\documentclass[sn-basic]{sn-jnl}% Math and Physical Sciences 
\RequirePackage{tikz}
\RequirePackage{fix-cm}
\usepackage{graphicx}
\usepackage{multirow}
\usepackage{amsthm}
\usepackage{amsmath,amssymb,amsfonts}
\usepackage{bbm}
\usepackage{mathrsfs}
\usepackage[title]{appendix}
\usepackage{xcolor}
\usepackage{textcomp}
\usepackage{manyfoot}
\usepackage{booktabs}
\usepackage{algorithm}
\usepackage{algorithmicx}
\usepackage{algpseudocode}
\usepackage{listings}
\usepackage[noabbrev]{cleveref}
\usepackage{siunitx}
\usepackage{graphicx}
\usepackage{xcolor}
\usepackage{pgfplots}
\usepackage{subfig}
\usepackage{stfloats}
\usepackage{caption}
\usepackage{soul}

\usetikzlibrary{matrix}
\pgfplotsset{compat=1.8}
\usepgfplotslibrary{groupplots}
\usepgfplotslibrary{colorbrewer}

\definecolor{seabornBlue}{RGB}{72, 120, 208}
\definecolor{seabornOrange}{RGB}{230, 159, 0}
\definecolor{seabornRed}{RGB}{214, 95, 95}
\definecolor{seabornGreen}{RGB}{106, 204, 100}
\definecolor{yellow}{RGB}{238, 204, 22}
\definecolor{teal}{RGB}{0, 129, 118}
\definecolor{darkRed}{rgb}{0.7569, 0.1529, 0.1765}
\definecolor{darkBlue}{rgb}{0, 0, 0.6549}

\DeclareMathOperator*{\argmin}{arg\,min}

\newtheorem{remark}{Remark}%

\theoremstyle{thmstylethree}%

\begin{document}

% \title[Article Title]{Control of water distribution networks with time-coupling constraints via distributed nonconvex optimization}

\title[Article Title]{Distributed nonconvex optimization for control of water networks with time-coupling constraints}

\author*[1]{\fnm{Bradley} \sur{Jenks}}\email{b.jenks21@imperial.ac.uk}

% \author[1]{\fnm{Aly-Joy} \sur{Ulusoy}}
% \author[2]{\fnm{Filippo} \sur{Pecci}}
% \author[1]{\fnm{Ivan} \sur{Stoianov}}

\author[1]{\fnm{Aly-Joy} \sur{Ulusoy}}\email{aly-joy.ulusoy15@imperial.ac.uk}
\author[2]{\fnm{Filippo} \sur{Pecci}}\email{fp0820@princeton.edu}
\author[1]{\fnm{Ivan} \sur{Stoianov}}\email{ivan.stoianov@imperial.ac.uk}

\affil*[1]{\orgdiv{Civil and Environmental Engineering}, \orgname{Imperial College London}, \city{London}, \postcode{SW7 2BB}, \country{United Kingdom}}

\affil[2]{\orgdiv{Andlinger Center for Energy and the Environment}, \orgname{Princeton University}, \state{New Jersey}, \postcode{08544}, \country{United States}}

\date{2023}

%%==================================%%
%%              ABSTRACT            %%
%%==================================%%

\abstract{In this paper, we present a new control model for optimizing pressure and water quality operations in water distribution networks. Our formulation imposes a set of time-coupling constraints to manage temporal pressure variations, which are exacerbated by the transition between pressure and water quality controls. The resulting optimization problem is a nonconvex, nonlinear program with nonseparable structure across time steps. This problem proves challenging for state-of-the-art nonlinear solvers, often precluding their direct use for near real-time control in large-scale networks. To overcome this computational burden, we investigate a distributed optimization approach based on the alternating direction method of multipliers (ADMM). In particular, we implement and evaluate two algorithms: a standard ADMM scheme and a two-level variant that provides theoretical convergence guarantees for our nonconvex problem. We use a benchmarking water network and a large-scale operational network in the UK for our numerical experiments. The results demonstrate good convergence behavior across all problem instances for the two-level algorithm, whereas the standard ADMM approach struggles to converge in some instances. With an appropriately tuned penalty parameter, however, both distributed algorithms yield good quality solutions and computational times compatible with near real-time (e.g. hourly) control requirements for large-scale water networks.}


\maketitle

%%==================================%%
%%                BODY              %%
%%==================================%%

%%% Introduction %%%
\section{Introduction}
\label{sec:introduction}

% Introduction to theme and general research problem %
Water distribution networks (WDNs) are undergoing a shift towards dynamic control, enabling them to adapt to various operational objectives in near real-time \citep{GIUDICIANNI2020,BUI2022,ULUSOY2023}. This shift is essential to maintain optimal operations amid the increasing challenges posed by asset deterioration, regulatory compliance, and constrained resources. A critical concern for water companies is the management of discolouration risk, which can significantly impact water quality and customer satisfaction \citep{BOXALL2023}. This concern has motivated a growing body of research on the design and control of self-cleaning networks for mitigating discolouration risk \citep{ABRAHAM2016,ABRAHAM2018,JENKS2023a}. Recently, a multi-objective optimization problem investigated the integration of self-cleaning operations with pressure management, a widely adopted practice for reducing leakage in WDNs \citep{JENKS2023b}. A major challenge identified was the considerable fluctuations in operating pressures that could arise due to conflicting hydraulic states associated with pressure and self-cleaning controls. These dynamics raise concerns for water companies, as they must weigh the benefits of dynamic control against the reliability challenges of aging infrastructure. 

% Literature review on engineering context to support pressure variation constraint %
There is growing evidence that suggests mean operating pressure and pressure variation influence pipe break rates. For example, \citet{REZAEI2015} found a positive correlation between pressure variation (up to $25 \, \si{\meter}$) and breaks caused by longitudinal cracks in operational WDNs in the UK. Similarly, \citet{MARTINEZ2016} identified pressure range as the most statistically significant factor in pipe breaks. \citet{ARRIAGADA2021} demonstrated a $61\%$ reduction in cast iron pipe breaks for a $10 \, \si{\meter}$ decrease in pressure range, highlighting the benefits of active pressure management. These findings are relevant to our work, as changes in hydraulic states resulting from pressure and self-cleaning controls have been shown to significantly affect operating pressures \citep{WRIGHT2015,ABRAHAM2016,JENKS2023b}. Given its role in infrastructure reliability, pressure variation has also been considered as an objective in the design of dynamically adaptive WDNs \citep{PECCI2017b}.

% Literaure review on decomposition-coordination methods %
Introducing pressure variation constraints in the control problem couples otherwise independent time steps. The resulting time-coupled, nonconvex optimization problem often precludes the direct application of state-of-the-art nonlinear solvers. Specifically, these solvers become incompatible with near real-time control implementations in large-scale networks which require fast and scalable solution methods. In the literature, similar optimization problems have been addressed in the context of pump scheduling where tank levels are coupled by flow at consecutive time steps. For example, \citet{ZESSLER1989} and \citet{NITIV1996} applied temporal decomposition and an iterative optimality technique to solve the time-coupled pump scheduling problem. In \citet{GHADDAR2015}, a Lagrange decomposition method was integrated with a simulation-based feasibility search. Similarly, the alternating direction method of multipliers (ADMM) was introduced in \citet{FOOLADIVANDA2018} to decouple feasibility constraints from the objective of minimizing pumping energy costs. ADMM was also applied in \citet{ZAMZAM2019} to decompose a coupled water-power network flow problem. 

ADMM is a well-known decomposition-coordination technique for solving large-scale convex optimization problems using distributed computing architectures \citep{BOYD2010}. Our problem is particularly suitable for ADMM since it exhibits a block structure that can be decomposed into smaller subproblems. These subproblems are efficiently distributed among different computing agents, while a coordination step ensures information consensus at each iteration. However, the presence of nonconvexity in many engineering applications, including the operation of water networks, may prevent a standard ADMM implementation from converging to a stationary solution \citep{WANG2019}. Recently, the literature has explored theoretical convergence properties of ADMM for different classes of nonconvex problems. The majority of these studies, however, concern problems where nonconvexity resides solely in the objective function \citep{LI2015,HONG2016,THEMELIS2018,WANG2019}. For nonconvex constraints, ADMM convergence has only been established under relatively strong assumptions \citep{MAGNUSSON2016}, which are often difficult to prove in practice. A promising approach to overcome these convergence issues, while preserving the desirable properties of ADMM, is the two-level distributed algorithm proposed in \citet{SUN2023}. This ADMM variant introduces slack variables in the consensus constraint to form a structure compliant with known ADMM convergence properties. An outer level augmented Lagrangian method (ALM) then drives the slack variables to zero using an amplifying penalty scheme. Recent work has successfully applied this two-level ADMM variant to optimal power flow problems \citep{SUN2021,GHOLAMI2023} and large-scale model predictive control (MPC) implementations \citep{TANG2022}. 
% In this paper, we implement both a standard ADMM scheme and the two-level algorithm proposed in \citet{SUN2023} for the control of water networks with time-coupling constraints.

% Paper contributions %
This paper presents a new control model for optimizing pressure and water quality operations in dynamically adaptive WDNs. Our formulation imposes a set of time-coupling constraints to manage pressure variations caused by the transition between pressure and water quality control objectives. We then investigate a distributed optimization approach to solve the resulting time-coupled, nonconvex problem, as state-of-the-art nonlinear solvers struggle to find feasible solutions. In particular, we implement a standard ADMM scheme and a two-level variant that provides convergence guarantees for our nonconvex problem. We demonstrate the computational performance of the distributed algorithms using a benchmarking water network and a large-scale operational network in the UK. Our analysis explores the empirical convergence behaviour of different ADMM penalty parameters, highlighting the robustness of the two-level algorithm to varying algorithmic parameters. With an appropriately tuned penalty parameter, however, both distributed algorithms exhibit superior computational performance compared to a centralized solution approach using a state-of-the-art non-linear optimization solver. We conclude discussing the application of distributed optimization for near real-time control in large-scale water networks. 

% Paper organization %
The rest of the paper is organized as follows. \Cref{sec:problem_formulation} introduces the problem formulation, including the network model, operational objectives, and time-coupling constraints. In \Cref{sec:distributed_algorithm}, we describe the standard ADMM and two-level distributed algorithms implemented to solve the time-coupled, nonconvex problem. In \Cref{sec:numerical_experiments}, we present case study data and the computational setup of numerical experiments. Lastly, in \Cref{sec:results}, we evaluate the performance of the distributed algorithms and discuss considerations for their practical implementation.

%%% Problem formulation %%%
\section{Problem formulation}
\label{sec:problem_formulation}
In the following, we formulate an optimization problem for scheduling pressure and self-cleaning controls in water networks.

%%% Water network model subsection %%%
\subsection{Water network model}   
\label{sec:water_network}
We model a water network as a directed graph $\mathcal{G}(\mathcal{N}, \mathcal{P})$. The set $\mathcal{N}$ comprises $n_n$ junction and $n_0$ source nodes and the set $\mathcal{P}$ comprises $n_p$ pipe or valve links. Network connectivity is represented through link-node incidence matrices $A_{12} \in \mathbb{R}^{n_p \times n_n}$ and $A_{10} \in \mathbb{R}^{n_p \times n_0}$ for junction and source nodes, respectively. We define $A_{12}$ (and $A_{10}$) using the following convention
\begin{equation} \label{eq:incidence_matrix}
    \begin{alignedat}{3}
        &A_{12}(j,i) =
        \begin{cases} 
            1 &\text{if link $j$ enters node $i$} \\
            0 & \text{if link $j$ is not connected to node $i$} \\
            -1 &\text{if link $j$ leaves node $i$.}
        \end{cases}
    \end{alignedat}
\end{equation}

We consider a steady-state hydraulic model to simulate hydraulic states. Let ${\mathcal{T} = \{1,\dots, n_t\}}$ denote a finite control horizon, which represents a collection of discrete time steps spanning a moving window of size $n_t$. For each time step $t \in \mathcal{T}$, the vectors of link flows ${q_t \in \mathbb{R}^{n_p}}$ and hydraulic heads ${h_t \in \mathbb{R}^{n_n}}$ denote unknown hydraulic states. Known hydraulic conditions are given by vectors of nodal demands ${d_t \in \mathbb{R}^{n_n}}$ and source hydraulic heads ${h_{0t} \in \mathbb{R}^{n_0}}$. Moreover, we introduce control via (i) pressure control valves (PCVs), which modulate pressure in a single flow direction, and (ii) automatic flushing valves (AFVs), which discard water at designated locations. The vectors ${\eta_t \in \mathbb{R}^{n_v}}$ and ${\alpha_t \in \mathbb{R}^{n_f}}$ model local losses across $n_v$ PCV links and operational demands at $n_f$ AFV nodes, respectively. Matrices ${A_{13} \in \mathbb{R}^{n_p \times n_v}}$ and ${A_{14} \in \mathbb{R}^{n_n \times n_f}}$ map known actuator locations. 

Hydraulic states $q_t$ and $h_t$ are governed by the following energy \eqref{eq:hyd_energy} and mass \eqref{eq:hyd_mass} conservation equations at time $t$:
\begin{subequations}
    \begin{align}
        &A_{12}h_t + A_{10}h_{0t} + \phi(q_t) + A_{13}\eta_t = 0,
        \label{eq:hyd_energy} \\
        &A_{12}^T q_t - d_t - A_{14}\alpha_t = 0, \label{eq:hyd_mass}
    \end{align}
\end{subequations}
where $\phi(\cdot)$ is a nonlinear function modelling frictional head loss across each link $j \in \mathcal{P}$, defined as 
\begin{equation}\label{eq:head_loss_model}
    \phi_j(q_{j,t}) = r_j|q_{j,t}|^{n_j-1}q_{j,t}.
\end{equation}
For valve links, $n_j = 2$ and 
\begin{equation} \label{eq:local_loss_model}
    r_j = \frac{8K_j}{g\pi^2D_j^4},
\end{equation}
where $K_j$ and $D_j$ denote the valve loss coefficient and diameter, respectively at link $j$. For pipe links, we apply the empirical Hazen-Williams formula, which has parameters $n_j = 1.852$ and 
\begin{equation}
    \label{eq:HW_model}
    r_j = \frac{10.67L_j}{C_j^{1.852}D_j^{4.871}},
\end{equation}
where $L_j$ is the length of link $j$; $C_j$ is the dimensionless H-W coefficient; and $D_j$ is the diameter of link $j$. Note that the semi-empirical Darcy-Weisbach formula can also be used to compute head loss across pipe links with no change to the problem formulation.

We introduce upper and lower bounds on continuous variables $h_t$, $\eta_t$, and $\alpha_t$ to define the feasible solution space for all $t \in \mathcal{T}$,
\begin{subequations}
\label{eq:variable_bounds}
\begin{align}
    &h^{\min}_t \leq h_t \leq h^{\max}_t, \label{eq:hyd_bounds_b}\\
    &\eta^L_t \leq \eta_t \leq \eta^U_t, \label{eq:hyd_bounds_c}\\
    &0 \leq \alpha_t \leq \alpha^{\max}_t, \label{eq:hyd_bounds_d}
\end{align}
\end{subequations}
where $h^{\min}_t$ is set to the minimum regulatory pressure and ${h^{\max}_t = \max_{i = 1, \dots, n_0}(h_{0t})_i \times \mathbbm{1}_{n_n}}$; $\eta^L_t$ and $\eta^U_t$ are established from corresponding $h^{\min}_t$ and $h^{\max}_t$ values at the upstream and downstream nodes of PCV links; and $\alpha^{\max}_t$ sets the maximum flushing rate at AFV nodes. Note that $h^{\min}_t$ and $\alpha^{\max}_t$ are defined on the basis of local regulatory and network conditions. Additionally, we enforce control direction at PCV links by introducing the following bi-linear constraint on flow $q_t$ and local loss $\eta_t$ variables for all $t \in \mathcal{T}$:
\begin{equation} \label{eq:pcv_direction}
    q_{j,t}(A_{13}\eta_t)_j \geq 0, \quad \forall j \in \mathcal{P}.
\end{equation}

%%% Control objectives subsection %%%
\subsection{Control objectives}
\label{sec:control_objectives}
In this work, we focus on the coordination of two control objectives. The first objective concerns hydraulic pressure, which plays a critical role in background leakage and pipe burst frequency \citep{SCHWALLER2015}. We define the network's average zone pressure (AZP) as a surrogate measure of pressure-induced leakage and pipe circumferential stress \citep{WRIGHT2015}, expressed as the linear function 
\begin{equation} \label{eq:azp_objective}
    f_{\text{AZP}}(h_t) := \sum_{i=1}^{n_n} w_{\text{AZP}i} (h_{i,t} - \zeta_i),
\end{equation}
where $h_{i,t}$ is the computed hydraulic head for node $i$ at time step $t$; $\zeta_i$ is the elevation of node $i$; and $w_{\text{AZP}i}$ is a coefficient weighting node $i$ by the length of its connected links \citep[Equation 5]{WRIGHT2015}.

The second objective aims to reduce discolouration risk via self-cleaning velocities. Recent research has investigated different optimization problems for maximizing the self-cleaning capacity (SCC) of a network \citep{ABRAHAM2016, ABRAHAM2018, JENKS2023a}. The SCC objective function is defined as the length of pipe with flow velocities exceeding a threshold sufficient to mobilize accumulated material \citep{VREEBURG2009,BLOKKER2010}. In order to use gradient-based optimization methods, \citet{ABRAHAM2016} proposed a smooth sum of logistic functions to approximate the SCC objective. This nonconvex approximation is written as  
\begin{equation} \label{eq:scc_objective}
    f_{\text{SCC}}(q_t) := \sum_{j=1}^{n_p}w_{\text{SCC}j} \left(\psi^+_{j}\left(\frac{q_{j,t}}{s_j}\right) + \psi^-_j\left(\frac{q_{j,t}}{s_j}\right)\right),
\end{equation}
where $\psi^+_{j}$ and $\psi^-_{j}$ are logistic functions in the positive and negative flow directions, respectively, for link $j$ (see \citet{JENKS2023a} for details); $q_{j,t}$ is the flow conveyed by link $j$ at time step $t$; $s_j$ is the cross-sectional area of link $j$; and $w_{\text{SCC}j}$ is a coefficient normalizing the length of link $j$ to the entire network. 

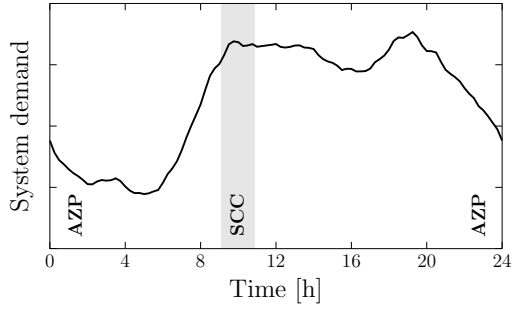
\begin{figure}[h]
    \centering
    \begin{tikzpicture}[scale=0.6]
    \begin{axis}[xlabel={Time [h]}, ylabel={System demand}, ymin={0}, ymax={80}, yticklabels={none}, xmin={0}, xmax={24}, xtick={{0, 4, ..., 24}}, tick style={black}, width={11.5cm}, height={7.0cm}, label style = {font=\Large}, tick label style = {font=\large},]
        \draw[draw={none}, fill=black, fill opacity={0.1}, label={SCC mode}] ({axis cs:9.1,0}|-{rel axis cs:0,1}) rectangle ({axis cs:10.9,0}|-{rel axis cs:0,0});
        \draw[style={ultra thick}, draw={none}, fill={none}, label={AZP mode}] ({axis cs:0.3,0}|-{rel axis cs:0,1}) rectangle ({axis cs:8.8,0}|-{rel axis cs:0,0});
        \draw[style={ultra thick}, draw={none}, fill={none}] ({axis cs:11.2,0}|-{rel axis cs:0,1}) rectangle ({axis cs:23.7,0}|-{rel axis cs:0,0});
        \addplot+[style={solid, very thick}, mark={none}, color={black}]
            coordinates {
                (0.0,35.26841712463647)
                (0.25,31.254816244356334)
                (0.5,28.829025129554793)
                (0.75,27.497198650613427)
                (1.0,25.915668422356248)
                (1.25,24.609983446542174)
                (1.5,23.65494578354992)
                (1.75,22.459373749792576)
                (2.0,21.056941949296743)
                (2.25,20.971345910569653)
                (2.5,21.890368898399174)
                (2.75,22.395459552295506)
                (3.0,22.18605152028612)
                (3.25,22.333163103554398)
                (3.5,22.939340411918238)
                (3.75,22.02742602257058)
                (4.0,20.282253063982353)
                (4.25,18.91351426835172)
                (4.5,18.119421026669443)
                (4.75,18.168997154571116)
                (5.0,17.802996183512732)
                (5.25,18.052574360743165)
                (5.5,18.663620290812105)
                (5.75,19.14579897443764)
                (6.0,21.487662029685453)
                (6.25,24.33414126257412)
                (6.5,26.178114238427952)
                (6.75,28.58496055356227)
                (7.0,32.134212164906785)
                (7.25,36.31122280983254)
                (7.5,40.052274737507105)
                (7.75,43.5312851651106)
                (8.0,47.07227364415303)
                (8.25,51.868349047610536)
                (8.5,56.414166253991425)
                (8.75,58.82891071261838)
                (9.0,60.29794540209696)
                (9.25,63.17634987598285)
                (9.5,66.66651481832378)
                (9.75,67.64093343680725)
                (10.0,67.37148238462396)
                (10.25,66.11641402356327)
                (10.5,66.29797395481728)
                (10.75,66.72271741344593)
                (11.0,65.87155898590572)
                (11.25,66.09951942088082)
                (11.5,66.29133809986524)
                (11.75,66.50045822327957)
                (12.0,66.75798537023365)
                (12.25,65.81241849204525)
                (12.5,65.61404157872312)
                (12.75,65.82745698327199)
                (13.0,66.27247556531802)
                (13.25,66.33457813691348)
                (13.5,65.62180143943988)
                (13.75,65.52399056451395)
                (14.0,65.11004951829091)
                (14.25,63.162039656192064)
                (14.5,61.84588004113175)
                (14.75,61.071451340103515)
                (15.0,60.62815871788189)
                (15.25,59.71080931369215)
                (15.5,58.273070556577295)
                (15.75,58.69758847076446)
                (16.0,58.63838244625367)
                (16.25,57.82668587937951)
                (16.5,57.83926808810793)
                (16.75,57.93165975622833)
                (17.0,58.74386182194576)
                (17.25,60.80642498959787)
                (17.5,61.831338291754946)
                (17.75,63.02237841812894)
                (18.0,65.00274277222343)
                (18.25,67.78994224080816)
                (18.5,68.91264707781374)
                (18.75,68.7620232028421)
                (19.0,69.6523414880503)
                (19.25,70.6682156403549)
                (19.5,69.04678659653291)
                (19.75,66.29048642446287)
                (20.0,64.47132633370347)
                (20.25,64.42739237402566)
                (20.5,63.999280128162354)
                (20.75,60.97791475849226)
                (21.0,58.35287849674933)
                (21.25,57.0642217125278)
                (21.5,55.95853077666834)
                (21.75,54.61441090912558)
                (22.0,52.42391668399796)
                (22.25,50.564290234819055)
                (22.5,49.13868002966046)
                (22.75,46.503933703526855)
                (23.0,45.150913272518665)
                (23.25,43.232921628979966)
                (23.5,41.05430295551196)
                (23.75,38.91701353690587)
                (24.0,35.26841712463647)
            }
            ;
        \node[above, black, rotate=90] at (20,110) {\large{\textbf{AZP}}};
        \node[above, black, rotate=90] at (106,110) {\large{\textbf{SCC}}};
        \node[above, black, rotate=90] at (234,110) {\large{\textbf{AZP}}};
    \end{axis}
    \end{tikzpicture}
    \vspace{0.25cm}
    \caption{AZP (primary) and SCC (secondary) control modes across the 24-hour control horizon $\mathcal{T}$. Note that the SCC window (shaded area) can be modified to align with network specific conditions.}
    \label{fig:control_horizon}
\end{figure}

Given that background leakage is a continuous phenomenon, we designate AZP as the primary (or continuous) control mode. Conversely, we set SCC as the secondary (or periodic) control mode, activated within a predefined 1-hour window. This coordination strategy is based on the assumption that achieving self-cleaning velocities daily (or every other day) mitigates the risk of discolouration in WDNs \citep{VREEBURG2009,BLOKKER2010}. \Cref{fig:control_horizon} illustrates the 24-hour control horizon, with system demand overlaid to highlight SCC activation during peak demand.

Here, we activate SCC controls between 09:30 to 10:30 to leverage the high flow velocities from peak demands. However, the duration and/or activation of this window can be easily modified in the problem formulation. Additionally, the SCC window may vary depending on the spatial distribution of demand characteristics in the network, particularly when targeting self-cleaning conditions for a specific subset of links.

%%% Time-coupling PV constraint subsection %%%
\subsection{Time-coupling constraints}
\label{sec:time-coupling_constraints}
We introduce a set of time-coupling constraints to manage nodal pressure range across the 24-hour control horizon. These constraints are defined by the following inequality:
\begin{equation}
\label{eq:pr_constraint}
    \max_{t\in \mathcal{T}}(h_{i,t}) - \min_{t\in \mathcal{T}}(h_{i,t}) \leq \delta, \quad \forall i \in \mathcal{N},
\end{equation}
where $\delta$ corresponds to a specified pressure range tolerance. For ease of notation, we introduce the set $\bar{\mathcal{X}} = \{x \in \mathbb{R}^{n_t \times n_n} \; | \; \max_{t\in \mathcal{T}}(x_{i,t}) - \min_{t\in \mathcal{T}}(x_{i,t}) \leq \delta, \; \forall i \in \mathcal{N}\}$. Note that, in our practical implementations, we reformulate \eqref{eq:pr_constraint} as a set of linear constraints by introducing auxiliary variables $c_{ui}$ and $c_{li}$ and setting
\begin{subequations}
\label{eq:pr_constraint_linear}
    \begin{alignat}{3}
    &h_{i,t} \leq c_{ui}, \quad \forall t \in \mathcal{T} \label{eq:pr_constraint_a}\\
    &h_{i,t} \geq c_{li}, \quad \forall t \in \mathcal{T} \label{eq:pr_constraint_b} \\
    &c_{ui} - c_{li} \leq \delta \label{eq:pr_constraint_c}.
    \end{alignat}
\end{subequations}

The choice of pressure range tolerance $\delta$ is based on its observed impact on pipe breaks, as demonstrated in statistical analyses from previous studies \citep{MARTINEZ2016, REZAEI2015, ARRIAGADA2021}. We explore the effect of $\delta$ on the resulting control actions ($u$) and hydraulic states ($q, h$) in \Cref{sec:results}. Additionally, we do not consider the frequency of pressure variation events since our problem focuses on changes in daily (steady-state) operating pressures. System dynamics were omitted on the assumption that control inputs in practice are designed to be gradual and thus mitigate the occurrence of sudden pressure changes.

%%% Overall problem formulation %%%
\subsection{Control problem}
\label{sec:overall_problem_formulation}
This work aims to optimize the AZP-SCC control schedule described in \Cref{sec:control_objectives} while managing temporal pressure variations through time-coupling constraints described in \Cref{sec:time-coupling_constraints}. The corresponding optimization problem can be formulated in compact form as
\begin{subequations}
\label{eq:control_problem}
\begin{alignat}{3}
    & \text{minimize}
    & \quad & \sum_{t \in \mathcal{T}} f_t(q_t, h_t) \label{eq:control_obj} \\
    & \text{subject to}
    & & (q_t, h_t, u_t) \in \mathcal{X}_t, \label{eq:control_hyd} \quad \forall t \in \mathcal{T} \\
    & & & h \in \bar{\mathcal{X}}, \label{eq:control_coupling} 
\end{alignat}
\end{subequations}
where $f_t := -f_{\text{SCC}}$ if $t \in \mathcal{T}_{\text{SCC}}$ and $f_t := f_{\text{AZP}}$ if $t \notin \mathcal{T}_{\text{SCC}}$; $u_t = (\eta_t, \alpha_t)$ defines the control variables at time $t$; $\mathcal{X}_t$ collects the hydraulic constraints described in \eqref{eq:hyd_energy}, \eqref{eq:hyd_mass}, \eqref{eq:variable_bounds}, and \eqref{eq:pcv_direction}; and $\bar{\mathcal{X}}$ represents the set of time-coupling constraints in \eqref{eq:pr_constraint}.

Problem \eqref{eq:control_problem} is a nonconvex, nonlinear programming (NLP) problem whose size grows with the number of discrete time steps in $\mathcal{T}$. A large number of time steps may be required to capture the diurnal periodicity of loading conditions as well as leverage the high-resolution frequencies of modern network sensing. In addition, time-coupling constraints \eqref{eq:control_coupling} introduce a nonseparable structure. These characteristics make problem \eqref{eq:control_problem} challenging to solve for near real-time control implementations in large-scale water networks.

%%% Distributed optimization algorithm %%%
\section{Distributed optimization}
\label{sec:distributed_algorithm}
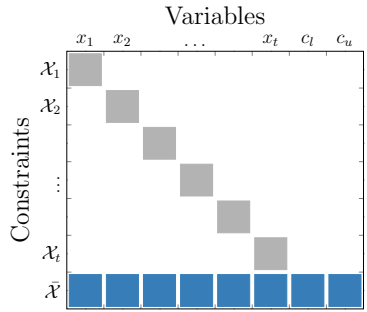
\begin{figure}[h]
    \centering
    \begin{tikzpicture}[scale=0.6]
    \begin{axis}[
        axis equal image,
        % scatter,
        xmin=1,
        xmax=8,
        xtick = {1, 2, ..., 8},
        xticklabels = {$x_1$, $x_2$, , $\dots$, , $x_t$, $c_l$, $c_u$},
        ymin=1,
        ymax=7,
        ytick = {1, 2, ..., 7},
        yticklabels = {$\mathcal{X}_1$, $\mathcal{X}_2$, , $\vdots$, , $\mathcal{X}_t$, $\bar{\mathcal{X}}$},
        minor tick num=1,
        tickwidth=0pt,
        y dir=reverse,
        xticklabel pos=right, % X axis labels go to the top
        enlargelimits={abs=0.5}, % Add half a unit on all sides
        xlabel = Variables,
        ylabel = Constraints,
        label style = {font=\Large},
        tick label style = {font=\large},
        cycle list/Set1-3,
    ]
    \addplot +[
        only marks,
        mark = square*,
        mark size = 10,
        index of colormap=1 of Set1-3,
        ] 
        coordinates {
            (1,7)
            (2,7)
            (3,7)
            (4,7)
            (5,7)
            (6,7)
            (7,7)
            (8,7)
        }
        ;
    \addplot +[
        only marks,
        mark = square*,
        mark size = 10,
        mark options={fill=black, fill opacity=0.3, draw=black, opacity=0.3},
        ] 
        coordinates {
            (1,1)
            (2,2)
            (3,3)
            (4,4)
            (5,5)
            (6,6)
        }
        ;
    \end{axis}
    \end{tikzpicture}
    
    \vspace{0.25cm}
    \caption{Block structure of problem \eqref{eq:control_problem} with time-coupling constraints $\bar{\mathcal{X}}$.}
    \label{fig:block_structure} 
\end{figure}

Problem \eqref{eq:control_problem} exhibits a special structure in which the only coupling between time steps is present in constraints \eqref{eq:control_coupling}. This nearly separable block structure, illustrated in Figure \ref{fig:block_structure}, is well-suited for decomposition-coordination procedures. In this work, we apply the alternating direction of multipliers (ADMM) algorithm, a widely adopted decomposition-coordination procedure for handling problems with coupling constraints \citep{ECKSTEIN1992,BOYD2010}.

We first describe a standard ADMM implementation by introducing a global variable copy (or variable splitting) in \Cref{sec:standard_admm}. This decomposes problem \eqref{eq:control_problem} into a series of smaller subproblems independent of time, enabling the use of distributed computing to enhance computational performance. Since the standard ADMM implementation does not provide theoretical convergence guarantees for nonconvex problems, we also implement the two-level ADMM variant \citep{SUN2021,SUN2023} in \Cref{sec:two_level_algorithm}, which is designed to handle problems with nonconvex constraints.

The following notation is used in the subsequent sections. The identity matrix of size $n$ is denoted by $I_n$. The inner product of $x, y \in \mathbb{R}^n$ is denoted by $x^Ty$ or $\langle x, y \rangle$. We use $\|x\| := \sqrt{\langle x, x \rangle}$ to denote the Euclidean norm of $x \in \mathbb{R}^n$. The indicator function of constraint set $\mathcal{X}$ is denoted by $\mathbb{I}_{\mathcal{X}}(x)$, with $\mathbb{I}_{\mathcal{X}}(x) = 0$ if $x \in \mathcal{X}$ and $+\infty$ if $x \notin \mathcal{X}$. The projection operator onto a closed set $\mathcal{X}$ is denoted as $\operatorname{Proj}_{\mathcal{X}}(x)$. Lastly, we define $[n] = \{1, \dots, n\}$.

\subsection{Standard ADMM}
\label{sec:standard_admm}
We first present a distributed reformulation of problem \eqref{eq:control_problem} by introducing a global copy of hydraulic state variables $\bar{h} \in \mathbb{R}^{n_n \times n_t}$ and the consensus constraint
\begin{equation}
\label{eq:global_copy}
   h - \bar{h} = 0,
\end{equation}
After duplicating variables, problem~\eqref{eq:control_problem} is equivalently written as
\begin{subequations}
\label{eq:problem_reform_1}
\begin{alignat}{3}
    & \text{minimize}
    & \quad & \sum_{t \in \mathcal{T}} f_t(q_t,h_t) \label{eq:prob_reform_1a} \\
    & \text{subject to}
    & & h - \bar{h} = 0, \label{eq:prob_reform_1b} \\
    & & & (q_t, h_t, u_t) \in \mathcal{X}_t, \label{eq:prob_reform_1c} \quad \forall t \in \mathcal{T} \\
    & & & \bar{h} \in \bar{\mathcal{X}}, \label{eq:prob_reform_1d}
\end{alignat}
\end{subequations}
which decomposes the problem across time steps as time-coupling constraints $\bar{\mathcal{X}}$ are now strictly a function of duplicated variables $\bar{h}$.

A standard ADMM implementation minimizes the augmented Lagrangian of problem \eqref{eq:problem_reform_1} \citep{BOYD2010}, defined as
\begin{equation} \label{eq:lagrangian_admm}
\begin{aligned}
        L(q,h,u,\bar{h},y) = &\sum_{t \in \mathcal{T}} \Big\{ f_t(q_t,h_t,u_t) + \mathbb{I}_{\mathcal{X}_t}(q_t,h_t,u_t) + \langle y_t, h_t - \bar{h}_t \rangle \\
        & + \frac{\rho}{2} \|h_t - \bar{h}_t\|^2 \Big\} + \mathbb{I}_{\bar{\mathcal{X}}}(\bar{h}) \\
        = & \sum_{t \in \mathcal{T}} L_t(q_t,h_t,u_t,\bar{h}_t,y_t) + \mathbb{I}_{\bar{\mathcal{X}}}(\bar{h}),
\end{aligned}
\end{equation}
where $\rho > 0$ is a fixed penalty parameter; $y_t \in \mathbb{R}^{n_n}$ is the vector of dual variables associated with consensus constraint \eqref{eq:prob_reform_1b}, with vectors $y = [\{y_t\}_{t \in \mathcal{T}}]^T$; and $\mathbb{I}_{\mathcal{X}_t}(\cdot)$ and $\mathbb{I}_{\bar{\mathcal{X}}}(\cdot)$ are the indicator functions of hydraulic feasibility set $\mathcal{X}_t$ and time-coupling constraint set $\bar{\mathcal{X}}$, respectively. Given initial values $(q^k,h^k,u^k,\bar{h}^k,y^k)$, ADMM minimizes \eqref{eq:lagrangian_admm} by performing a sequence of steps over $k$ iterations, as follows:
\begin{subequations}
\label{eq:admm_algorithm}
    \begin{alignat}{3}
    q_t^{k+1},h_t^{k+1}, u_t^{k+1} &:= \underset{\substack{q_t,h_t,u_t}}{\argmin} \; L_t(q_t,h_t,u_t,\bar{h}_t^k,y_t^k), \,\; \forall t \in \mathcal{T} \label{eq:admm_x_update} \\
    \bar{h}^{k+1} &:= \underset{\substack{\bar{h}}}{\argmin} \; L(q^{k+1},h^{k+1},u^{k+1},\bar{h},y^k) \label{eq:admm_hbar_update} \\
    y^{k+1} &:= y^k + \rho(h^{k+1} - \bar{h}^{k+1}). \label{eq:admm_dual_update}
    \end{alignat}
\end{subequations}

The sequence of ADMM steps can be efficiently handled by state-of-the-art optimization solvers. With $\bar{h}^k$ and $y^k$ fixed, locally optimal state $(q_t^{k+1},h_t^{k+1})$ and control $u_t^{k+1}$ variables can be independently computed for each nonconvex NLP subproblem in \eqref{eq:admm_x_update}. The NLP subproblems can therefore be distributed across multiple processes to improve computational performance (see \Cref{sec:numerical_experiments} for distributed computing implementation used in this work). In the subsequent step, $h^{k+1}$ is broadcast to the coordination problem in \eqref{eq:admm_hbar_update}, which is a relatively small convex quadratic problem enforcing time-coupling constraints $\bar{h} \in \bar{\mathcal{X}}$. Lastly, dual variables $y^{k+1}$ are updated in \eqref{eq:admm_dual_update}, with the fixed penalty parameter $\rho$ as the step size.

We terminate ADMM if solution $(q^k,h^k,u^k,\bar{h}^k,y^k)$ at iteration $k$ satisfies the following criteria:
\begin{subequations}
    \begin{alignat}{3}
    &\|\rho (\bar{h}^{k+1} - \bar{h}^k)\| \leq \epsilon_{\text{d}} \label{eq:dual_residual} \\
    &\|h^k - \bar{h}^k\| \leq \epsilon_{\text{p}}, \label{eq:primal_residual}
    \end{alignat}
\end{subequations}
where \eqref{eq:dual_residual} and \eqref{eq:primal_residual} denote the dual feasibility and primal residuals, respectively; and $\epsilon_{\text{d}}$ and $\epsilon_{\text{p}}$ are positive tolerances.  

Recent literature has explored theoretical convergence properties for nonconvex ADMM problems \citep{HONG2016,MAGNUSSON2016,WANG2019,JIANG2019}. However, \cite{SUN2021,SUN2023} demonstrate that problem \eqref{eq:problem_reform_1} fails to satisfy the sufficient conditions necessary for ensuring convergence to a stationary solution. Specifically, the local constraints on the last block update in \eqref{eq:admm_hbar_update} prevent the use of the unconstrained optimality condition to link primal and dual variables. Although these conditions guarantee convergence, they are not necessarily required for ADMM to reach a stationary solution. Thus, we empirically examine the convergence behavior of the standard ADMM scheme from a practical perspective in \Cref{sec:results}.

\subsection{Two-level algorithm}
\label{sec:two_level_algorithm}
Here, we implement the two-level ADMM variant proposed in \citet{SUN2023}. Unlike the standard approach, the two-level algorithm provides theoretical convergence guarantees to a stationary solution under mild assumptions for the nonconvex setting.

We consider the following reformulation of problem \eqref{eq:problem_reform_1} with slack variables $z \in \mathbb{R}^{n_n \times n_t}$,
\begin{subequations}
\label{eq:problem_reform_2}
\begin{alignat}{3}
    & \text{minimize}
    & \quad & \sum_{t \in \mathcal{T}} f_t(q_t,h_t) \label{eq:prob_reform_2a} \\
    & \text{subject to}
    & & h - \bar{h} + z = 0 \label{eq:prob_reform_2b} \\
    & & & (q_t, h_t, u_t) \in \mathcal{X}_t, \label{eq:prob_reform_2c} \quad \forall t \in \mathcal{T} \\
    & & & \bar{h} \in \bar{\mathcal{X}} \label{eq:prob_reform_2d} \\
    & & & z = 0. \label{eq:prob_reform_2d}
\end{alignat}
\end{subequations}
In this reformulation, the consensus constraint \eqref{eq:prob_reform_2b} involves three variable blocks $(h, \bar{h}, z)$. To preserve ADMM convergence properties, \citet{SUN2023} consider constraint $z=0$ separately through an augmented Lagrangian relaxation (ALR) of problem \eqref{eq:problem_reform_2}. The ALR at iteration $m$ is formulated as
\begin{subequations}
\label{eq:problem_alr}
\begin{alignat}{3}
    & \text{minimize}
    & \quad & \sum_{t \in \mathcal{T}} \Big\{ f_t(q_t,h_t) + \langle \lambda_t^m, z_t \rangle + \frac{\beta^m}{2} \|z_t\|^2 \Big\} \label{eq:prob_alr_a} \\
    & \text{subject to}
    & & h - \bar{h} + z = 0 \label{eq:prob_alr_b} \\
    & & & (q_t, h_t, u_t) \in \mathcal{X}_t, \label{eq:prob_alr_c} \quad \forall t \in \mathcal{T} \\
    & & & \bar{h} \in \bar{\mathcal{X}}, \label{eq:prob_alr_d}
\end{alignat}
\end{subequations}
where $\lambda_t^m \in \mathbb{R}^{n_n}$ is the vector of dual variables associated with constraint \eqref{eq:prob_reform_2d} at iteration $m$, with vectors $\lambda = [\{\lambda_t\}_{t \in \mathcal{T}}]^T$; and $\frac{\beta^m}{2} \|z_t\|^2$ is a quadratic penalty term with $\beta^m >0$ at iteration $m$ and time step $t$. Relaxing constraint $z=0$, we can solve problem \eqref{eq:problem_alr} via a three-block ADMM scheme that now satisfies the unconstrained optimality condition of the last variable block $z$.

Given $\lambda^m$ and $\beta^m$ at the $m$-th outer level iteration, ADMM minimizes the augmented Lagrangian function associated with problem \eqref{eq:problem_alr}, defined as 
\begin{equation} \label{eq:lagrangian_two_level}
\begin{aligned}
    \hat{L}^m(q,h,u,\bar{h},z,y) = &\sum_{t \in \mathcal{T}} \Big\{ f_t(q_t,h_t,u_t) + \mathbb{I}_{\mathcal{X}_t}(q_t,h_t,u_t) + \langle \lambda_t^m, z_t \rangle \\
    & + \frac{\beta^m}{2} \|z_t\|^2 + \langle y_t, h_t - \bar{h}_t + z_t \rangle \\
    & + \frac{\rho^m}{2} \|h_t - \bar{h}_t + z_t\|^2 \Big\} + \mathbb{I}_{\bar{\mathcal{X}}}(\bar{h}) \\
    = & \sum_{t \in \mathcal{T}} \hat{L}^m_t(q_t,h_t,u_t,\bar{h}_t,z_t,y_t) + \mathbb{I}_{\bar{\mathcal{X}}}(\bar{h}),
\end{aligned}
\end{equation} 
where $y_t \in \mathbb{R}^{n_n}$ is the vector of dual variables associated with constraint \eqref{eq:prob_alr_b}, with vectors $y = [\{y_t\}_{t \in \mathcal{T}}]^T$; and $\rho^m>0$ is the ADMM penalty parameter at iteration $m$ of the outer level problem. The three-block ADMM updates $(q^k,h^k,u^k,\bar{h}^k,z^k,y^k)$ in sequence over $k$ inner level iterations, with $\lambda^m$, $\beta^m$, and $\rho^m$ kept constant at each outer level iteration $m$:
\begin{subequations}
\label{eq:3block_admm_algorithm}
    \begin{alignat}{3}
    q_t^{k+1},h_t^{k+1}, u_t^{k+1} &:= \underset{\substack{q_t,h_t,u_t}}{\argmin} \; \hat{L}_t(q_t,h_t,u_t,\bar{h}_t^k,z_t^k,y_t^k), \,\; \forall t \in \mathcal{T} \label{eq:3block_admm_x_update} \\
    \bar{h}^{k+1} &:= \underset{\substack{\bar{h}}}{\argmin} \; \hat{L}(q^{k+1},h^{k+1},u^{k+1},\bar{h},z^k,y^k) \label{eq:3block_admm_hbar_update} \\
    z^{k+1} &:= \underset{\substack{z}}{\argmin} \; \hat{L}(q^{k+1},h^{k+1},u^{k+1},\bar{h}^{k+1},z,y^k) \label{eq:3block_admm_z_update} \\
    y^{k+1} &:= y^k + \rho^m(h^{k+1} - \bar{h}^{k+1} + z^{k+1}). \label{eq:3block_admm_dual_update}
    \end{alignat}
\end{subequations}
Similar to the standard ADMM implementation in \eqref{eq:admm_algorithm}, hydraulic state $(q_t^{k+1},h_t^{k+1})$ and control $u_t^{k+1}$ variables are decoupled across time steps $t \in \mathcal{T}$ in \eqref{eq:3block_admm_x_update} and can thus be updated in a distributed manner. Under mild assumptions, which we show problem \eqref{eq:problem_reform_2} satisfy in Appendix~\ref{sec:A1}, \citet{SUN2023} prove that the inner level ADMM described in \eqref{eq:3block_admm_algorithm} converges to an approximate stationary solution of problem \eqref{eq:problem_alr} when $\rho^m = 2\beta^m$.

In practice, we terminate ADMM if solution $(q^k,h^k,u^k,\bar{h}^k,z^k,y^k)$ at inner level iteration $k$ satisfies the following stopping criteria:
\begin{subequations}
\label{eq:admm_termination}
    \begin{alignat}{3}
    &\|\rho^m (\bar{h}^{k-1} + z^{k-1} - \bar{h}^k - z^k)\| \leq \epsilon_1 \label{eq:inner_residual_2} \\
    &\|\rho^m (z^{k-1} - z^k)\| \leq \epsilon_2 \label{eq:inner_residual_2} \\
    &\|h^k - \bar{h}^k + z^k\| \leq \epsilon_3. \label{eq:inner_residual_3}
    \end{alignat}
\end{subequations}
where $\epsilon_{i \in [3]}$ are positive tolerances.

However, constraint \eqref{eq:prob_reform_2d} may not necessarily be satisfied after ADMM successfully terminates at the inner level. A classical augmented Lagrangian method is therefore adopted to drive slack variables $z$ to zero by updating $\lambda^{m+1}$, as follows:
\begin{equation}
\label{eq:alm_update}
   \lambda^{m+1} = \operatorname{Proj}_{[\underline{\lambda}, \overline{\lambda}]}(\lambda^m + \beta^m z^m),
\end{equation}
in which the projection onto a predetermined hypercube $[\underline{\lambda}, \overline{\lambda}]$ is introduced to ensure boundedness of dual variables and thereby the augmented Lagrangian function $\hat{L}$. We also amplify the penalty parameter $\beta^m$ if $z^m$ has not sufficiently decreased from the previous iteration $z^{m-1}$. The $\beta^{m+1}$ update step is written as
\begin{equation} \label{eq:beta_update}
    \begin{alignedat}{3}
        &\beta^{m+1} =
        \begin{cases} 
            \gamma \beta^m, \quad &\|z^m\| > \omega\|z^{m-1}\| \\
            \beta^m, \quad &\|z^m\| \leq \omega\|z^{m-1}\|
        \end{cases}
    \end{alignedat}
\end{equation}
where $\gamma>1$ and $\omega \in [0,1)$ are tunable parameters. We terminate the outer level ALM if solution $(q^m,h^m,u^m,\bar{h}^m,z^m)$ satisfies the primal residual stopping criterion
\begin{equation}
\label{eq:alm_termination}
   \|h^k - \bar{h}^k\| \leq \epsilon_p.
\end{equation}

In summary, we solve problem \eqref{eq:problem_reform_2} in two levels: the inner level applies a three-block ADMM to find an approximate stationary solution to problem \eqref{eq:problem_alr}, with iterates indexed by $k$; and the outer level drives slack variables $z$ to zero using a classic ALM framework, with iterates indexed by $m$. This two-level algorithm is repeated until the outer level stopping criterion \eqref{eq:alm_termination} is met. Pseudocode for the two-level algorithm is presented in \Cref{alg:two_level_algorithm}.
\begin{algorithm}[t]
  \caption{Two-level distributed algorithm}
  \label{alg:two_level_algorithm}
  \begin{algorithmic}[1]
    \State \textbf{Set} dual variable bounds $[\underline{\lambda},\overline{\lambda}]$; $\gamma > 1$; $\omega \in [0,1)$; and terminating tolerances $\epsilon_i > 0, \; \forall i \in [3]$
    \State \textbf{Initialize} starting points $(q_t^1, h_t^1, u_t^1) \in \mathcal{X}_t$ for all $t \in \mathcal{T}$ and $\bar{h}^1 \in \bar{\mathcal{X}}$; penalty parameter $\beta^1 > 0$; and dual variables $\lambda^1 \in [\underline{\lambda},\overline{\lambda}]$
    \State $m \leftarrow 1$
    \While{outer level stopping criterion \eqref{eq:alm_termination} is not met}
        \State \textbf{Initialize} starting points $(q^0, h^0, u^0) := (q^m, h^m, u^m)$, $\bar{h}^0 := h^m$, $z^0 := 0_{n_n \times n_t}$, and $y^0 := -\lambda^m$, such that $(q_t^0, h_t^0, u_t^0) \in \mathcal{X}_t$ for all $t \in \mathcal{T}$ and $\lambda^m + \beta^mz^0 + y^0 = 0$; and penalty parameter $\rho^m = 2\beta^m$ 
        \State $k \leftarrow 0$
        \While{inner level stopping criteria \eqref{eq:admm_termination} are not met}
            \State First block update (in parallel): $(q_t^{k+1}, h_t^{k+1},u_t^{k+1}), \; \forall t \in \mathcal{T}$ $\leftarrow$ step \eqref{eq:3block_admm_x_update}
            \State Second block update: $\bar{h}^{k+1}$ $\leftarrow$ step \eqref{eq:3block_admm_hbar_update}
            \State Third block update: $z^{k+1}$ $\leftarrow$ step \eqref{eq:3block_admm_z_update} 
            \State Inner dual update: $y^{k+1}$ $\leftarrow$ step \eqref{eq:3block_admm_dual_update}
            \State $k \leftarrow k + 1$
        \EndWhile
        \State $(q^m, h^m, u^m, \bar{h}^m, z^m, y^m) \leftarrow (q^k, h^k, u^m, \bar{h}^k, z^k, y^
        k)$
        \State Outer dual update: $\lambda^{m+1}$ $\leftarrow$ step \eqref{eq:alm_update}
        \State Penalty parameter update: $\beta^{m+1}$ $\leftarrow$ step \eqref{eq:beta_update}
        \State $m \leftarrow m + 1$
    \EndWhile
    \State Return $\epsilon$-stationary primal solution $(q^m,h^m,u^m,\bar{h}^m,z^m)$
  \end{algorithmic}
\end{algorithm}

%%% Numerical experiments %%%
\section{Numerical experiments}
\label{sec:numerical_experiments}

\subsection{Case studies and problem data}
\label{sec:case_networks}
\begin{figure*}[h]
    \centering
    \subfloat[\centering \label{fig:case_networks_modena}\texttt{Modena}]{{\includegraphics[scale=0.48]{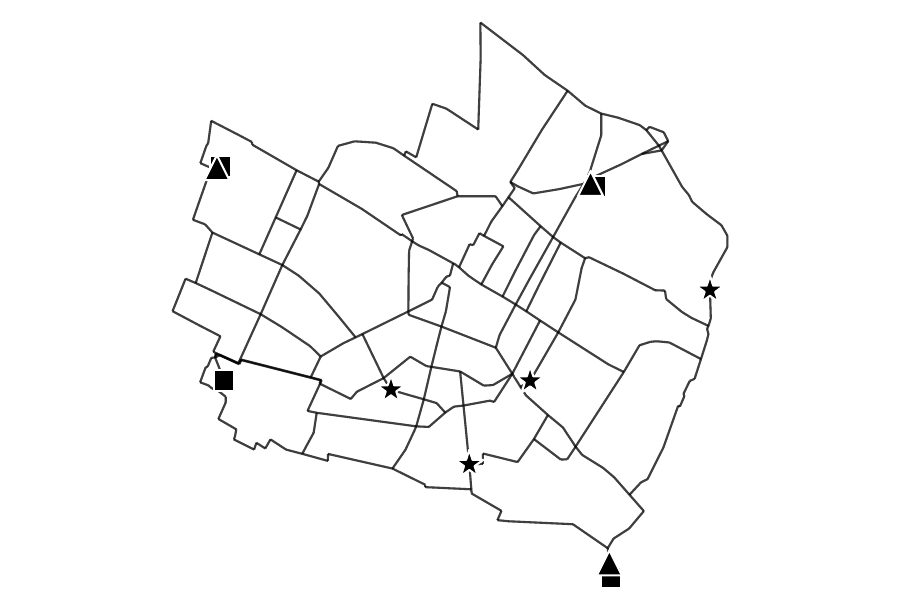} }} \\
    \vspace{0.25cm}
    \subfloat[\centering \label{fig:case_networks_bwfl}\texttt{BWFLnet}]{{\includegraphics[scale=0.56]{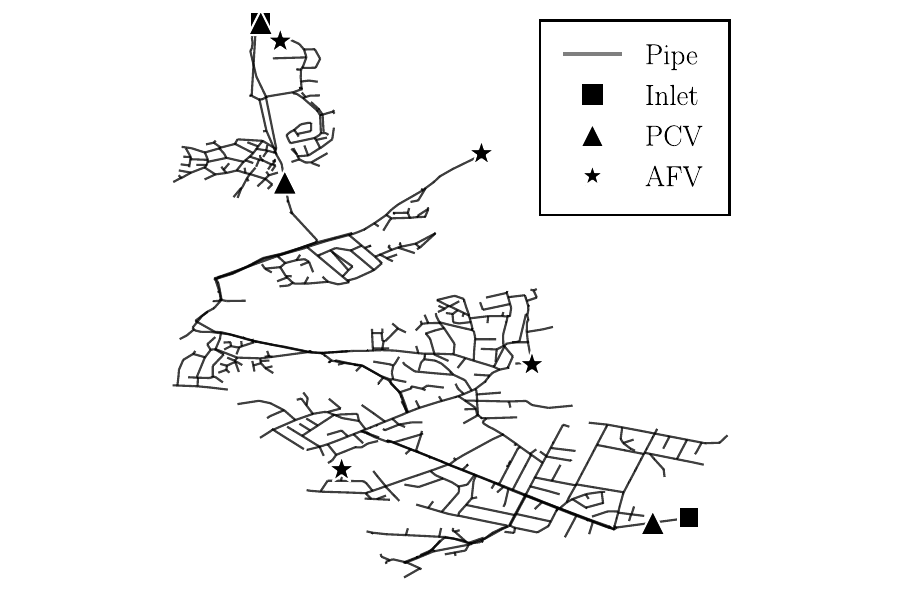} }}
    \vspace{0.25cm}
    \caption{Case network layouts.}
    \label{fig:case_networks}
\end{figure*}

We consider two case study networks: (i) \texttt{Modena}, a medium-sized benchmark model comprised of $n_0 = 4$ source nodes, ${n_n = 268}$ junction nodes, and ${n_p = 317}$ links \citep{BRAGALLI2012}; and (ii) Bristol Water Field Lab network (\texttt{BWFLnet}), a large-scale operational network model comprised of ${n_0 = 2}$ source nodes, ${n_n = 2745}$ junction nodes, and ${n_p = 2816}$ links \citep{WRIGHT2014}. Both models feature $n_v = 3$ PCV and $n_f = 4$ AFV control actuators, with locations chosen \textit{a priori} using the optimization method proposed in ~\citet{JENKS2023b}. \Cref{fig:case_networks} illustrates the networks' topology and control valve locations.

Furthermore, we consider a 24-hour control horizon $\mathcal{T}$ to capture the diurnal periodicity of loading conditions (i.e. demands) in WDNs. For \texttt{Modena}, we use $n_t=24$ 1-hour time steps, which is a typical resolution for hydraulic modelling and existing control implementations. On the other hand, we use $n_t=96$ 15-minute time steps in \texttt{BWFLnet} to leverage higher resolution sensor data increasingly available in operational networks. Accordingly, control implementations in practice will need to handle increasingly larger optimization problem sizes. The sizes of the resulting optimization problems investigated in this work are summarized in \Cref{tab:problem_size}.
\begin{table}[h]
\centering
\caption{Problem sizes for case study networks.}
\vspace{0.25cm}
\setlength{\tabcolsep}{10pt}
{\renewcommand{\arraystretch}{1.2}
\small{
\begin{tabular}{lccc}
\toprule
Network & \# Cont. var. & \# Lin. cons. & \# Noncon. terms \\ \midrule
\texttt{Modena}    & $14,452$   & $47,932$    & $8,314$ \\ [3pt]
\texttt{BWFLnet}   & $537,369$   & $1,862,553$     & $284,800$ \\ [3pt]
\bottomrule
\end{tabular}
}
}
\label{tab:problem_size}
\end{table}

The following parameters were chosen to formulate the time-coupled, nonconvex problem \eqref{eq:control_problem}. We assume a 1-hour self-cleaning mode is activated during the morning peak demand period, as illustrated previously in \Cref{fig:control_horizon}. This is denoted by sets $\mathcal{T}_{\text{SCC}} = \{7, 8\}$ and $\mathcal{T}_{\text{SCC}} = \{38, \dots, 42\}$ for \texttt{Modena} and \texttt{BWFLnet}, respectively, where $\mathcal{T}_{\text{SCC}} \subset \mathcal{T}$. We applied a maximum flow rate $\alpha^{\max} = 25 \, \si{\liter/\second}$ at AFV actuators and a minimum regulatory pressure head of $h^{\min} = 15 \, \si{\meter}$ (UK regulations) at model nodes. Moreover, we considered three pressure range tolerances $\delta \in \{20, 15, 10\}$ (in meters) to define time-coupling constraint set $\bar{\mathcal{X}}$ in \eqref{eq:pr_constraint}. These tolerances were selected with reference to the unconstrained case (i.e. $\bar{\mathcal{X}}=\mathbb{R}$) and the pressure range results reported in \citet{REZAEI2015} and \citet{ARRIAGADA2021}.

\subsection{Implementation details}
\label{sec:implementation}
All numerical experiments were written in Julia 1.9.1 \citep{BEZANSON2017} with optimization solvers accessed via JuMP 1.12.0 \citep{DUNNING2017}. The code was compiled on a 2.50-GHz Intel(R) Core(TM) i9-11900H CPU with 8 cores and 32.0 GB of RAM running Ubuntu 20.04.5 LTS Linux distribution. Linear programs were solved using Gurobi 10.0.2 \citep{GUROBI2023}, and nonlinear programs were solved using IPOPT 3.14.4 \citep{WACHTER2006} with the HSL MA57 linear solver \citep{HSL2021}. Moreover, we used Julia's \texttt{Distributed.jl} module to distribute the computations in \eqref{eq:admm_x_update} and \eqref{eq:3block_admm_x_update} across eight computing agents, each running their own instance of Julia. This approach was found to be faster than multithreading on a single process.

The two-level algorithm described in \Cref{alg:two_level_algorithm} was implemented as follows. The dual variables $\lambda$ were bounded between $\underline{\lambda}=-10^{5}$ and $\overline{\lambda}=10^{5}$. We set $\gamma=1.25$ and $\omega=0.75$ for the penalty parameter $\beta^{m+1}$ update in \eqref{eq:beta_update}. To prevent scaling issues in IPOPT, we enforced a maximum penalty parameter of $\beta \leq 10^5$. The inner level ADMM penalty parameter $\rho$ was initialized as $\rho^m=2\beta^m$ at each outer level iteration $m$ \citep{SUN2023}, except for inner level iteration $k=0$ where $\rho^0 = 0$. We considered initial $\beta^1$ values ranging from $10^{-3}$ to $10^2$ (see results in \Cref{sec:results}). Both the standard and two-level algorithms were initialized with feasible starting points $(q_t^1, h_t^1, u_t^1) \in \mathcal{X}_t$, for all $t \in \mathcal{T}$. We generated this solution from hydraulic states corresponding to the no control case (i.e. $u_t^1 = 0$, for all $t \in \mathcal{T}$). Moreover, the duplicated variables were initialized to $\bar{h}^1 = h^1$, satisfying the time-coupling constraints $\bar{\mathcal{X}}$ since pressure range tolerances $\delta$ were chosen on the basis of initial hydraulic conditions. For each outer level iteration $m$, ADMM was reinitialized using the previous (feasible) solution $(q^0, h^0, u^0) := (q^{m-1}, h^{m-1}, u^{m-1})$, and by setting $\bar{h}^0 = h^{m-1}$, $z^0=0$, and $y^0=-\lambda^{m-1}$. We terminate the inner level ADMM when ${\|h^k - \bar{h}^k + z^k\| \leq \sqrt{n_n n_t} \mathbin{/} (100 \cdot m)}$ or ${\|\rho^m (z^{k-1} - z^k)\| \leq 10^{-5}}$. The former criterion ensures that we satisfy primal feasibility, with a stopping tolerance that becomes increasingly strict with the outer iteration count $m$, while the latter checks the stability of slack variable $z^k$ to inform early ADMM termination. The standard ADMM and outer level ALM of the two-level method terminate when primal feasibility satisfies $\|h^m - \bar{h}^m\| \leq \sqrt{n_n n_t} \cdot \epsilon_p$, where $\epsilon_p=10^{-2}$. This tolerance value was selected to align with the precision of pressure monitoring devices and the inherent uncertainties within hydraulic models.

\section{Results and discussion}
\label{sec:results}

\subsection{Centralized solver}
\label{sec:centralized_results}
We first discuss the performance of state-of-the-art nonlinear solver IPOPT \citep{WACHTER2006} for solving \eqref{eq:control_problem} in a centralized manner. The results for different pressure tolerances $\delta$ are summarized in \Cref{tab:results_comparison}. For the smaller \texttt{Modena} network, IPOPT successfully converged to a feasible solution in all problem instances with a maximum computational time on the order of $100$ seconds. On the other hand, IPOPT failed to compute feasible solutions within the 1-hour time limit for the larger \texttt{BWFLnet} network when time-coupling constraints $\bar{\mathcal{X}}$ were imposed. This outcome was expected due to the substantial size of the resulting NLP problem formulated for \texttt{BWFLnet} with $n_t=96$ discrete time steps (see \Cref{tab:problem_size}). Specifically, the dense nature of time-coupling constraints $\bar{\mathcal{X}}$ results in large memory requirements, making it challenging for IPOPT to handle efficiently. We also tested control periods with a smaller number of time steps for \texttt{BWFLnet} to better understand IPOPT's limitations. Interestingly, there was a significant jump in computational time from $345$ seconds to $> 3600$ seconds between problem instances with $n_t=24$ and $n_t=32$ time steps, respectively.

These results reveal that a centralized solution approach may not be practical for large-scale water networks. Consequently, distributed optimization becomes appealing to facilitate computationally efficient control strategies. Our experiments indicate that the computational bottleneck of centralized solvers is mainly due to the linear system solver's memory usage within the interior-point solver. In this respect, first-order methods (e.g. sequential convex programming~\citep{WRIGHT2015}) or matrix-free solvers (e.g. ALPAQA \citep{PAS2022}) may offer viable alternatives to efficiently solve problem \eqref{eq:control_problem}.

\subsection{Distributed optimization}
\label{sec:distributed_results}
\begin{figure}[p]
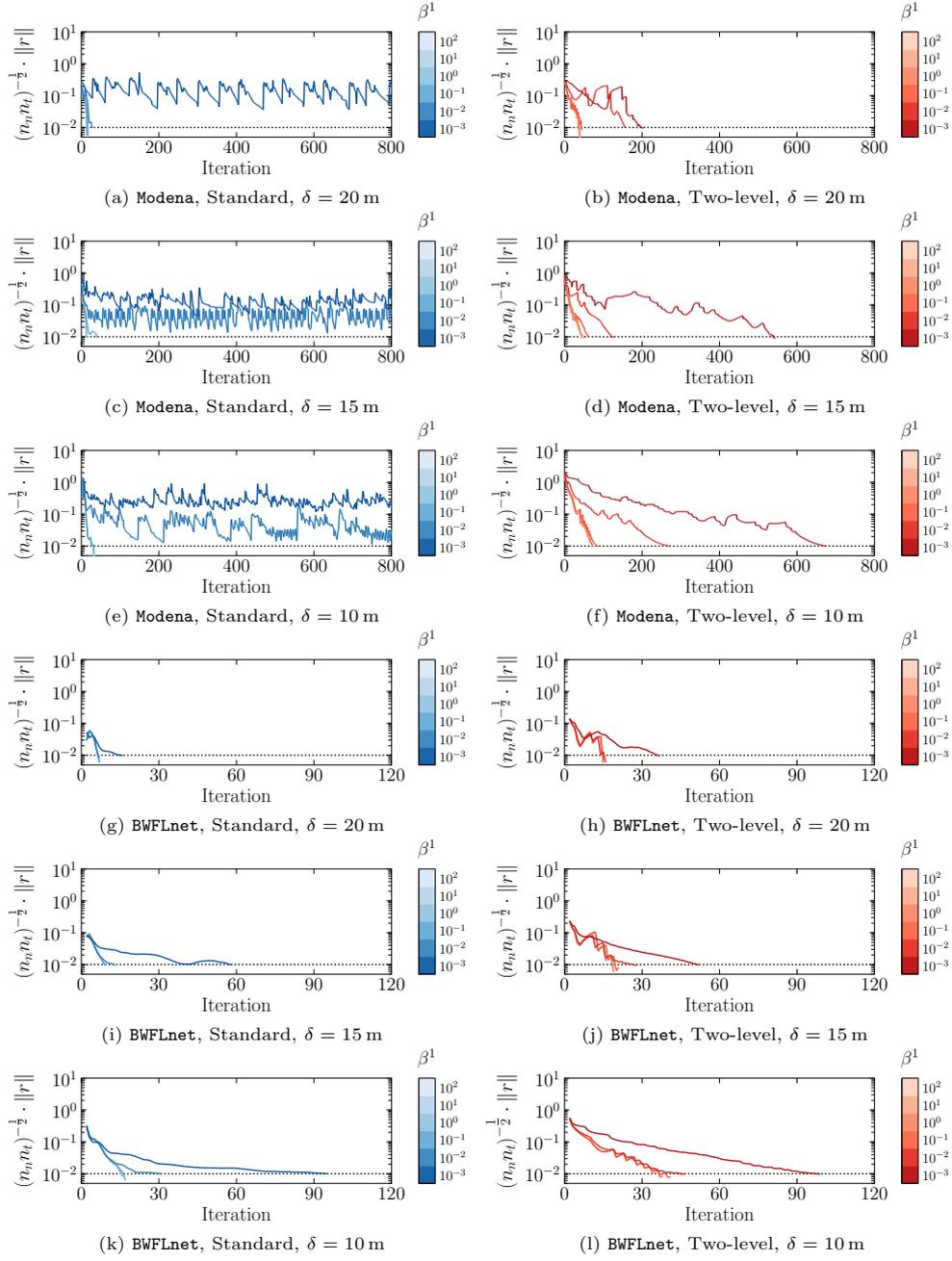

    \centering
    \subfloat[\centering \label{fig:residual_a}\texttt{Modena}, Standard, $\delta=20 \, \si{\meter}$]{\input{modena_residuals_admm_a}}
    \hfill
    \subfloat[\centering \label{fig:residual_b}\texttt{Modena}, Two-level, $\delta=20 \, \si{\meter}$]{\input{modena_residuals_two_level_a}}
    \\ \vspace{-0.4cm}
    \subfloat[\centering \label{fig:residual_c}\texttt{Modena}, Standard, $\delta=15 \, \si{\meter}$]{\input{modena_residuals_admm_b}}
    \hfill
    \subfloat[\centering \label{fig:residual_d}\texttt{Modena}, Two-level, $\delta=15 \, \si{\meter}$]{\input{modena_residuals_two_level_b}}
    \\ \vspace{-0.4cm}
    \subfloat[\centering \label{fig:residual_e}\texttt{Modena}, Standard, $\delta=10 \, \si{\meter}$]{\input{modena_residuals_admm_c}}
    \hfill
    \subfloat[\centering \label{fig:residual_f}\texttt{Modena}, Two-level, $\delta=10 \, \si{\meter}$]{\input{modena_residuals_two_level_c}}
    \\ \vspace{-0.4cm}
    \subfloat[\centering \label{fig:residual_g}\texttt{BWFLnet}, Standard, $\delta=20 \, \si{\meter}$]{% \pgfplotsset{
%     colormap={myblues}{
%         rgb255(0cm)=(255,255,255);
%         rgb255(1cm)=(222,235,247);
%         rgb255(2cm)=(198,219,239);
%         rgb255(3cm)=(158,202,225);
%         rgb255(4cm)=(107,174,214);
%         rgb255(5cm)=(49,130,189);
%         rgb255(6cm)=(8,81,156);
%     },
% }

\begin{tikzpicture}[scale=0.45]
\begin{axis}[
    width=11cm,
    height=4.8cm,
    ylabel={$(n_n n_t)^{-\frac{1}{2}} \cdot \|r\|$},
    xlabel={Iteration},
    xmin={0},
    xmax={120}, 
    xtick={{0, 30, ..., 120}},
    ymode={log},
    ymin={0.005}, 
    ymax={10}, 
    tick style={black},
    tick label style={{font=\Large}},
    label style={{font=\LARGE}},
    colormap/Blues-6, 
    colormap={reverse Blues-6}{
        indices of colormap={
            \pgfplotscolormaplastindexof{Blues-6},...,0 of Blues-6}
    },
    point meta min=-3,
    point meta max=3,
    colorbar right,
    colorbar style={
        title={$\beta^1$},
        title style={font={\Large}},
        tick label style={{font=\large}},
        ytick={-2.5, -1.5, -0.5, 0.5, 1.5, 2.5},
        yticklabels={$10^{-3}$,$10^{-2}$,$10^{-1}$,$10^{0}$,$10^{1}$,$10^2$},
        samples=7,
        tick style={draw=none}
    },
    colorbar sampled,
    ]
    \draw[style={dotted, very thick}, color={black}] ({rel axis cs:1,0}|-{axis cs:0,0.010000000000000002}) -- ({rel axis cs:0,0}|-{axis cs:0,0.010000000000000002})
    ;
    \addplot+[style={solid, very thick}, mesh, mark={none}, point meta={log10(1e2)}]
        coordinates {

            (2,0.030050850083481667)
            (3,0.06088833705318804)
            (4,0.051259662285843284)
            (5,0.029374497276468892)
            (6,0.012709836935054215)
            (7,0.005858980810627762)
            }
            ;
    \addplot+[style={solid, very thick}, mesh, mark={none}, point meta={log10(1e1)}]
        coordinates {

            (2,0.03005000534045795)
            (3,0.06088476733832881)
            (4,0.05125520163650384)
            (5,0.02938036152214277)
            (6,0.012712929109353427)
            (7,0.005860621575722514)
            }
            ;
    \addplot+[style={solid, very thick}, mesh, mark={none}, point meta={log10(1e0)}]
        coordinates {

            (2,0.03004237074876654)
            (3,0.06084959687283283)
            (4,0.05121459359622529)
            (5,0.029436996264346937)
            (6,0.01274559206066655)
            (7,0.005878148158947002)
            }
            ;
    \addplot+[style={solid, very thick}, mesh, mark={none}, point meta={log10(1e-1)}]
        coordinates {

            (2,0.02997257153726606)
            (3,0.0605043761426189)
            (4,0.05082191195038498)
            (5,0.029999298514204232)
            (6,0.013071479279525094)
            (7,0.006051844481112234)
            }
            ;
    \addplot+[style={solid, very thick}, mesh, mark={none}, point meta={log10(1e-2)}]
        coordinates {

            (2,0.02968313724014822)
            (3,0.05674884175638351)
            (4,0.048112148712577456)
            (5,0.034566986562312316)
            (6,0.013041475705122518)
            (7,0.008321814495431451)
            }
            ;
    \addplot+[style={solid, very thick}, mesh, mark={none}, point meta={log10(1e-3)}]
        coordinates {

            (2,0.05292500738102283)
            (3,0.03507970155673388)
            (4,0.04135448112743249)
            (5,0.036647925341857505)
            (6,0.025391395997925784)
            (7,0.017974171395218865)
            (8,0.014240692009359656)
            (9,0.0133962971980947)
            (10,0.013266231097077311)
            (11,0.012780364205370457)
            (12,0.01206869587182139)
            (13,0.011152427761673661)
            (14,0.010340960279500128)
            (15,0.010042198249872422)
            (16,0.009800214301926822)
            }
            ;
    \end{axis}
\end{tikzpicture}}
    \hfill
    \subfloat[\centering \label{fig:residual_h}\texttt{BWFLnet}, Two-level, $\delta=20 \, \si{\meter}$]{% \pgfplotsset{
%     colormap={myblues}{
%         rgb255(0cm)=(255,255,255);
%         rgb255(1cm)=(222,235,247);
%         rgb255(2cm)=(198,219,239);
%         rgb255(3cm)=(158,202,225);
%         rgb255(4cm)=(107,174,214);
%         rgb255(5cm)=(49,130,189);
%         rgb255(6cm)=(8,81,156);
%     },
% }

\begin{tikzpicture}[scale=0.45]
\begin{axis}[
    width=11cm,
    height=4.8cm,
    ylabel={$(n_n n_t)^{-\frac{1}{2}} \cdot \|r\|$},
    xlabel={Iteration},
    xmin={0},
    xmax={120}, 
    xtick={{0, 30, ..., 120}},
    ymode={log},
    ymin={0.005}, 
    ymax={10}, 
    tick style={black},
    tick label style={{font=\Large}},
    label style={{font=\LARGE}},
    colormap/Reds-6, 
    colormap={reverse Reds-6}{
        indices of colormap={
            \pgfplotscolormaplastindexof{Reds-6},...,0 of Reds-6}
    },
    point meta min=0,
    point meta max=6,
    point meta min=-3,
    point meta max=3,
    colorbar right,
    colorbar style={
        title={$\beta^1$},
        title style={font={\Large}},
        tick label style={{font=\large}},
        ytick={-2.5, -1.5, -0.5, 0.5, 1.5, 2.5},
        yticklabels={$10^{-3}$,$10^{-2}$,$10^{-1}$,$10^{0}$,$10^{1}$,$10^2$},
        samples=7,
        tick style={draw=none}
    },
    colorbar sampled,
    ]
    \draw[style={dotted, very thick}, color={black}] ({rel axis cs:1,0}|-{axis cs:0,0.010000000000000002}) -- ({rel axis cs:0,0}|-{axis cs:0,0.010000000000000002})
    ;
    \addplot+[style={solid, very thick}, mesh, mark={none}, point meta={log10(1e2)}]
 coordinates {
            (2,0.14086068547166122)
            (3,0.09223097888103983)
            (4,0.06952072465378333)
            (5,0.029456072941165208)
            (6,0.01846745227890473)
            (7,0.02648645032350841)
            (8,0.037177067226193246)
            (9,0.0464630631444919)
            (10,0.0541176847360756)
            (11,0.023677593255164386)
            (12,0.03111109159425031)
            (13,0.03601065827486027)
            (14,0.0393206731187684)
            (15,0.006630022710038926)
        }
        ;
    \addplot+[style={solid, very thick}, mesh, mark={none}, point meta={log10(1e1)}]
        coordinates {
            (2,0.14086068540640673)
            (3,0.09223298322459555)
            (4,0.06952454677120191)
            (5,0.029460334009060837)
            (6,0.018468555605077264)
            (7,0.026483751494753762)
            (8,0.037172586085759)
            (9,0.046457598417197765)
            (10,0.0541113996259922)
            (11,0.023681191755848766)
            (12,0.031114999223557115)
            (13,0.03601516858819146)
            (14,0.03932487044409041)
            (15,0.006633245836438674)
        }
        ;
    \addplot+[style={solid, very thick}, mesh, mark={none}, point meta={log10(1e0)}]
coordinates {
            (2,0.14086068726503337)
            (3,0.09225271900237166)
            (4,0.0695626418372963)
            (5,0.029502972853028917)
            (6,0.01847722360812408)
            (7,0.026456900394440055)
            (8,0.03712731942713722)
            (9,0.04640248465975528)
            (10,0.05404767940321149)
            (11,0.023717287357990037)
            (12,0.03115424327576925)
            (13,0.036060410880312685)
            (14,0.039366755385033377)
            (15,0.006666015476271288)
        }
        ;
    \addplot+[style={solid, very thick}, mesh, mark={none}, point meta={log10(1e-1)}]
coordinates {
            (2,0.14086068630334744)
            (3,0.0924510305839653)
            (4,0.06994466326969502)
            (5,0.02993173194288506)
            (6,0.018575087860213373)
            (7,0.026192852168572728)
            (8,0.03668146313740884)
            (9,0.04586041221520156)
            (10,0.05342104079238292)
            (11,0.024088746007517266)
            (12,0.03155383025270266)
            (13,0.0365009346589833)
            (14,0.03965066290792194)
            (15,0.005754959296262271)
   
        }
        ;
    \addplot+[style={solid, very thick}, mesh, mark={none}, point meta={log10(1e-2)}]
coordinates {
            (2,0.1408606867916751)
            (3,0.0944031480205661)
            (4,0.07366953658572938)
            (5,0.034144937106261035)
            (6,0.020222157849867962)
            (7,0.024186435033850658)
            (8,0.03286125071696666)
            (9,0.04098581004831967)
            (10,0.04782336173955805)
            (11,0.027754248617691657)
            (12,0.03478312587467212)
            (13,0.038564877880043515)
            (14,0.010434674180217669)
            (15,0.012935336434521012)
            (16,0.005901555995022887)
        }
        ;
    \addplot+[style={solid, very thick}, mesh, mark={none}, point meta={log10(1e-3)}]
coordinates {
            (2,0.14086068761906415)
            (3,0.11079019387618669)
            (4,0.10130362590298653)
            (5,0.0675225290641529)
            (6,0.05374934999171737)
            (7,0.03588149086602021)
            (8,0.03407791823954611)
            (9,0.03763706221722837)
            (10,0.04243602591972318)
            (11,0.04725286253907979)
            (12,0.0516970123439814)
            (13,0.055484457613257925)
            (14,0.04747522497766126)
            (15,0.04458292068493277)
            (16,0.04271864545416571)
            (17,0.035043590489969606)
            (18,0.030459222696930636)
            (19,0.026948962840457107)
            (20,0.022892185643107774)
            (21,0.019491581890895882)
            (22,0.017775882963927885)
            (23,0.017198735575894766)
            (24,0.01728986106419099)
            (25,0.017741907978837856)
            (26,0.01811538358171538)
            (27,0.017544384266826182)
            (28,0.01754361317086108)
            (29,0.017468648714436402)
            (30,0.015953434465479942)
            (31,0.015264169639650263)
            (32,0.014337723782136545)
            (33,0.01276350942111719)
            (34,0.012007557313800075)
            (35,0.010823242666828447)
            (36,0.0099771810811203)
            (37,0.009601099923143659)
        }
        ;
    \end{axis}
\end{tikzpicture}}
    \\ \vspace{-0.4cm}
    \subfloat[\centering \label{fig:residual_i}\texttt{BWFLnet}, Standard, $\delta=15 \, \si{\meter}$]{% \pgfplotsset{
%     colormap={myblues}{
%         rgb255(0cm)=(255,255,255);
%         rgb255(1cm)=(222,235,247);
%         rgb255(2cm)=(198,219,239);
%         rgb255(3cm)=(158,202,225);
%         rgb255(4cm)=(107,174,214);
%         rgb255(5cm)=(49,130,189);
%         rgb255(6cm)=(8,81,156);
%     },
% }

\begin{tikzpicture}[scale=0.45]
\begin{axis}[
    width=11cm,
    height=4.8cm,
    ylabel={$(n_n n_t)^{-\frac{1}{2}} \cdot \|r\|$},
    xlabel={Iteration},
    xmin={0},
    xmax={120}, 
    xtick={{0, 30, ..., 120}},
    ymode={log},
    ymin={0.005}, 
    ymax={10}, 
    tick style={black},
    tick label style={{font=\Large}},
    label style={{font=\LARGE}},
    colormap/Blues-6, 
    colormap={reverse Blues-6}{
        indices of colormap={
            \pgfplotscolormaplastindexof{Blues-6},...,0 of Blues-6}
    },
    point meta min=-3,
    point meta max=3,
    colorbar right,
    colorbar style={
        title={$\beta^1$},
        title style={font={\Large}},
        tick label style={{font=\large}},
        ytick={-2.5, -1.5, -0.5, 0.5, 1.5, 2.5},
        yticklabels={$10^{-3}$,$10^{-2}$,$10^{-1}$,$10^{0}$,$10^{1}$,$10^2$},
        samples=7,
        tick style={draw=none}
    },
    colorbar sampled,
    ]
    \draw[style={dotted, very thick}, color={black}] ({rel axis cs:1,0}|-{axis cs:0,0.010000000000000002}) -- ({rel axis cs:0,0}|-{axis cs:0,0.010000000000000002})
    ;
    \addplot+[style={solid, very thick}, mesh, mark={none}, point meta={log10(1e2)}]
        coordinates {

            (2,0.07649256677169565)
            (3,0.09410559334784412)
            (4,0.06734599677190341)
            (5,0.036572824859307046)
            (6,0.026565052588662045)
            (7,0.019120618464273258)
            (8,0.013327132080832556)
            (9,0.009933759348243457)
        }
        ;
    \addplot+[style={solid, very thick}, mesh, mark={none}, point meta={log10(1e1)}]
        coordinates {

            (2,0.07648787278655411)
            (3,0.09408943367909832)
            (4,0.06733632366090347)
            (5,0.03657900786282035)
            (6,0.026562425374736846)
            (7,0.01911457979797659)
            (8,0.013327428810921043)
            (9,0.009932502647464947)
        }
        ;
    \addplot+[style={solid, very thick}, mesh, mark={none}, point meta={log10(1e0)}]
        coordinates {

            (2,0.07644153108966832)
            (3,0.09392891588539849)
            (4,0.06724141106659082)
            (5,0.0366436981540373)
            (6,0.026533580947994562)
            (7,0.01907565936822616)
            (8,0.013338751481702864)
            (9,0.009828930346981226)
        }
        ;
    \addplot+[style={solid, very thick}, mesh, mark={none}, point meta={log10(1e-1)}]
        coordinates {

            (2,0.0759827187630206)
            (3,0.09228862470160906)
            (4,0.06598496275778239)
            (5,0.0374443934445298)
            (6,0.0263827530796607)
            (7,0.019128162179662095)
            (8,0.013674528534386077)
            (9,0.009734383967778005)
        }
        ;
    \addplot+[style={solid, very thick}, mesh, mark={none}, point meta={log10(1e-2)}]
        coordinates {

            (2,0.07319686126626224)
            (3,0.08041050707451737)
            (4,0.05667390323538312)
            (5,0.0398710558121984)
            (6,0.02462139311834192)
            (7,0.018331001564203864)
            (8,0.01579086227055092)
            (9,0.013104716174281971)
            (10,0.011392541571308768)
            (11,0.010897301691437334)
            (12,0.0105599220842391)
            (13,0.009448148471999027)
        }
        ;
    \addplot+[style={solid, very thick}, mesh, mark={none}, point meta={log10(1e-3)}]
coordinates {

            (2,0.08308145079249274)
            (3,0.07005648245323773)
            (4,0.06011123351895494)
            (5,0.04978785844999559)
            (6,0.03890434946530899)
            (7,0.03446698111424016)
            (8,0.032340138621040196)
            (9,0.03128906839035403)
            (10,0.03053945251276487)
            (11,0.02966007697124462)
            (12,0.028907364336739967)
            (13,0.028696691719013588)
            (14,0.026936008463473056)
            (15,0.024484319961606046)
            (16,0.023947281324371406)
            (17,0.023552311335395548)
            (18,0.02289763817108882)
            (19,0.02170736230340827)
            (20,0.021412441497683874)
            (21,0.021112970022493474)
            (22,0.02101471384905036)
            (23,0.021157176981377764)
            (24,0.021239324298567105)
            (25,0.02114238404188285)
            (26,0.02088385938554683)
            (27,0.020617467815299596)
            (28,0.02024162999651913)
            (29,0.019677258657125603)
            (30,0.018836410033692325)
            (31,0.01763676332856418)
            (32,0.016376113276825133)
            (33,0.015359432993322859)
            (34,0.014320526130857715)
            (35,0.013294314799557986)
            (36,0.012348742672313864)
            (37,0.011532684877698441)
            (38,0.010870764044185913)
            (39,0.010481720483227005)
            (40,0.01029433214925463)
            (41,0.010082799548774583)
            (42,0.010320155650390591)
            (43,0.010695001043814253)
            (44,0.011264494945823258)
            (45,0.011924450051091271)
            (46,0.012425646203673537)
            (47,0.012778108936289306)
            (48,0.013038127291983035)
            (49,0.013128802152156166)
            (50,0.013334005889673034)
            (51,0.013029748362507932)
            (52,0.01279161981931696)
            (53,0.012026846504103761)
            (54,0.011871937521454356)
            (55,0.01155407901182844)
            (56,0.011006927723674297)
            (57,0.010334776325606247)
            (58,0.009782700012651545)
        }
        ;
    \end{axis}
\end{tikzpicture}}
    \hfill
    \subfloat[\centering \label{fig:residual_j}\texttt{BWFLnet}, Two-level, $\delta=15 \, \si{\meter}$]{% \pgfplotsset{
%     colormap={myblues}{
%         rgb255(0cm)=(255,255,255);
%         rgb255(1cm)=(222,235,247);
%         rgb255(2cm)=(198,219,239);
%         rgb255(3cm)=(158,202,225);
%         rgb255(4cm)=(107,174,214);
%         rgb255(5cm)=(49,130,189);
%         rgb255(6cm)=(8,81,156);
%     },
% }

\begin{tikzpicture}[scale=0.45]
\begin{axis}[
    width=11cm,
    height=4.8cm,
    ylabel={$(n_n n_t)^{-\frac{1}{2}} \cdot \|r\|$},
    xlabel={Iteration},
    xmin={0},
    xmax={120}, 
    xtick={{0, 30, ..., 120}},
    ymode={log},
    ymin={0.005}, 
    ymax={10}, 
    tick style={black},
    tick label style={{font=\Large}},
    label style={{font=\LARGE}},
    colormap/Reds-6, 
    colormap={reverse Reds-6}{
        indices of colormap={
            \pgfplotscolormaplastindexof{Reds-6},...,0 of Reds-6}
    },
    point meta min=-3,
    point meta max=3,
    colorbar right,
    colorbar style={
        title={$\beta^1$},
        title style={font={\Large}},
        tick label style={{font=\large}},
        ytick={-2.5, -1.5, -0.5, 0.5, 1.5, 2.5},
        yticklabels={$10^{-3}$,$10^{-2}$,$10^{-1}$,$10^{0}$,$10^{1}$,$10^2$},
        samples=7,
        tick style={draw=none}
    },
    colorbar sampled,
    ]
    \draw[style={dotted, very thick}, color={black}] ({rel axis cs:1,0}|-{axis cs:0,0.010000000000000002}) -- ({rel axis cs:0,0}|-{axis cs:0,0.010000000000000002})
    ;
    \addplot+[style={solid, very thick}, mesh, mark={none}, point meta={log10(1e2)}]
        coordinates {
            (2,0.2324720988347559)
            (3,0.14309264960249865)
            (4,0.10477681040005897)
            (5,0.047330798305416304)
            (6,0.03979080114304891)
            (7,0.05309357087065235)
            (8,0.06800121716232313)
            (9,0.08070144244914576)
            (10,0.09105251560332704)
            (11,0.09947760680416405)
            (12,0.10638746673724901)
            (13,0.026590542415230414)
            (14,0.03332630812184216)
            (15,0.037809377355283245)
            (16,0.013294833617068017)
            (17,0.01709738996789134)
            (18,0.019367163505090577)
            (19,0.0059067865547787896)
        }
        ;
    \addplot+[style={solid, very thick}, mesh, mark={none}, point meta={log10(1e1)}]
        coordinates {
            (2,0.23247209886381845)
            (3,0.14309576181540373)
            (4,0.10478306726901601)
            (5,0.04733709790629489)
            (6,0.03979443045077654)
            (7,0.05308975037846091)
            (8,0.06798700080193118)
            (9,0.08068824740388776)
            (10,0.09103357144527899)
            (11,0.09945349598168934)
            (12,0.10635915109342664)
            (13,0.02659627516860264)
            (14,0.03333552431794166)
            (15,0.037817635676884315)
            (16,0.013309807767461622)
            (17,0.01711457078345168)
            (18,0.019380756787536854)
            (19,0.005906893297357789)
        }
        ;
    \addplot+[style={solid, very thick}, mesh, mark={none}, point meta={log10(1e0)}]
coordinates {
            (2,0.23247210617030342)
            (3,0.14312647981127197)
            (4,0.10484509471790254)
            (5,0.04740182052074398)
            (6,0.039833859587700596)
            (7,0.05305317975246263)
            (8,0.06787044666621964)
            (9,0.08055576112320724)
            (10,0.09084502269020617)
            (11,0.09921419128424258)
            (12,0.10607571862483109)
            (13,0.026655330038254023)
            (14,0.03342833602646307)
            (15,0.03790040964007327)
            (16,0.013462231404416616)
            (17,0.017283812054286817)
            (18,0.019506891368260032)
            (19,0.005915804947710018)
        }
        ;
    \addplot+[style={solid, very thick}, mesh, mark={none}, point meta={log10(1e-1)}]
coordinates {
            (2,0.23247210364350293)
            (3,0.1434388529114229)
            (4,0.10546998723488511)
            (5,0.04805196427655171)
            (6,0.040168144545938665)
            (7,0.05268038569178911)
            (8,0.06695472092973906)
            (9,0.07912743300974738)
            (10,0.08894225966850222)
            (11,0.09683233214854717)
            (12,0.10323780981931677)
            (13,0.02728943419606861)
            (14,0.03395081863289473)
            (15,0.03838960695971206)
            (16,0.04133927708200158)
            (17,0.013936522587994086)
            (18,0.017045091905382243)
            (19,0.01851055638880022)
            (20,0.006430795615430422)
            (21,0.0081484408817906)
        }
        ;
    \addplot+[style={solid, very thick}, mesh, mark={none}, point meta={log10(1e-2)}]
coordinates {
            (2,0.23247210468256263)
            (3,0.14645184574874243)
            (4,0.11138841987637817)
            (5,0.054415602008385106)
            (6,0.042948650857657096)
            (7,0.051009319833664246)
            (8,0.06187078024943006)
            (9,0.07075122719947774)
            (10,0.07733288996399704)
            (11,0.08204522354905378)
            (12,0.033187786886501223)
            (13,0.04172125814178803)
            (14,0.045295479573562905)
            (15,0.04686055627347008)
            (16,0.016280449881058846)
            (17,0.018556479796311764)
            (18,0.01976927888118004)
            (19,0.01288798834252014)
            (20,0.013400222213177968)
            (21,0.012066382021832876)
            (22,0.012363200550271198)
            (23,0.011087010993656964)
            (24,0.011229338831564076)
            (25,0.010311949816527829)
            (26,0.010199511658000163)
            (27,0.00965001971512728)
            (28,0.009471130012864752)
        }
        ;
    \addplot+[style={solid, very thick}, mesh, mark={none}, point meta={log10(1e-3)}]
coordinates {
            (2,0.23247209971519825)
            (3,0.17045789246129664)
            (4,0.1553129235183422)
            (5,0.101350973982656)
            (6,0.08394460029152007)
            (7,0.07697807521548312)
            (8,0.07460418437567415)
            (9,0.0743737941445521)
            (10,0.07501724953100378)
            (11,0.06978272666722277)
            (12,0.06396310819490474)
            (13,0.06189898263149505)
            (14,0.05680463038922071)
            (15,0.05190996589761379)
            (16,0.049920234891620646)
            (17,0.045897644931845935)
            (18,0.04210618199559824)
            (19,0.04035618570074385)
            (20,0.03829968319328778)
            (21,0.03582315645471104)
            (22,0.03459501054523354)
            (23,0.03298107756939932)
            (24,0.03128720196470822)
            (25,0.030388155220934337)
            (26,0.02863324325817075)
            (27,0.027284577711160787)
            (28,0.026516661157962768)
            (29,0.025280218203311507)
            (30,0.02421568223595153)
            (31,0.02355114754459727)
            (32,0.022700688191071332)
            (33,0.021893530515970074)
            (34,0.021354621204086113)
            (35,0.020404442143273385)
            (36,0.019725487336178394)
            (37,0.019244798475716737)
            (38,0.01832702912323483)
            (39,0.017642897558574493)
            (40,0.017131662993431107)
            (41,0.0163602081747513)
            (42,0.015736494678473364)
            (43,0.015245588495862388)
            (44,0.014369433717082153)
            (45,0.013800004641271277)
            (46,0.013310280962275304)
            (47,0.012436121661465713)
            (48,0.011843504832431339)
            (49,0.01135204013815409)
            (50,0.010728174181282272)
            (51,0.010202701374831903)
            (52,0.009751025829034532)
        }
        ;
    \end{axis}
\end{tikzpicture}}
    \\ \vspace{-0.4cm}
    \subfloat[\centering \label{fig:residual_k}\texttt{BWFLnet}, Standard, $\delta=10 \, \si{\meter}$]{% \pgfplotsset{
%     colormap={myblues}{
%         rgb255(0cm)=(255,255,255);
%         rgb255(1cm)=(222,235,247);
%         rgb255(2cm)=(198,219,239);
%         rgb255(3cm)=(158,202,225);
%         rgb255(4cm)=(107,174,214);
%         rgb255(5cm)=(49,130,189);
%         rgb255(6cm)=(8,81,156);
%     },
% }

\begin{tikzpicture}[scale=0.45]
\begin{axis}[
    width=11cm,
    height=4.8cm,
    ylabel={$(n_n n_t)^{-\frac{1}{2}} \cdot \|r\|$},
    xlabel={Iteration},
    xmin={0},
    xmax={120}, 
    xtick={{0, 30, ..., 120}},
    ymode={log},
    ymin={0.005}, 
    ymax={10}, 
    tick style={black},
    tick label style={{font=\Large}},
    label style={{font=\LARGE}},
    colormap/Blues-6, 
    colormap={reverse Blues-6}{
        indices of colormap={
            \pgfplotscolormaplastindexof{Blues-6},...,0 of Blues-6}
    },
    point meta min=0,
    point meta min=-3,
    point meta max=3,
    colorbar right,
    colorbar style={
        title={$\beta^1$},
        title style={font={\Large}},
        tick label style={{font=\large}},
        ytick={-2.5, -1.5, -0.5, 0.5, 1.5, 2.5},
        yticklabels={$10^{-3}$,$10^{-2}$,$10^{-1}$,$10^{0}$,$10^{1}$,$10^2$},
        samples=7,
        tick style={draw=none}
    },
    colorbar sampled,
    ]
    \draw[style={dotted, very thick}, color={black}] ({rel axis cs:1,0}|-{axis cs:0,0.010000000000000002}) -- ({rel axis cs:0,0}|-{axis cs:0,0.010000000000000002})
    ;
    \addplot+[style={solid, very thick}, mesh, mark={none}, point meta={log10(1e2)}]
        coordinates {

            (2,0.3288557209507901)
            (3,0.14694092343292747)
            (4,0.13032292303914103)
            (5,0.12567986498721287)
            (6,0.10415832542021636)
            (7,0.0801578922461254)
            (8,0.05916872883477715)
            (9,0.04129733024073268)
            (10,0.03538192471375184)
            (11,0.02712294510700637)
            (12,0.023454064817403628)
            (13,0.020842127207216905)
            (14,0.01619944329973915)
            (15,0.011276799591584858)
            (16,0.008553907770134503)
        }
        ;
    \addplot+[style={solid, very thick}, mesh, mark={none}, point meta={log10(1e1)}]
        coordinates {

            (2,0.3288541140284185)
            (3,0.14693546745923625)
            (4,0.13030814405182697)
            (5,0.12567247440583224)
            (6,0.10415646478400482)
            (7,0.08017043703345691)
            (8,0.05916891341216104)
            (9,0.04129604783821757)
            (10,0.03537490491963888)
            (11,0.027112866446241324)
            (12,0.023428638879443745)
            (13,0.020846229047185336)
            (14,0.016227905550355082)
            (15,0.011301442968746182)
            (16,0.008530905913909754)
       }
        ;
    \addplot+[style={solid, very thick}, mesh, mark={none}, point meta={log10(1e0)}]
        coordinates {

            (2,0.32884294429522454)
            (3,0.14688180141619683)
            (4,0.13015732900202368)
            (5,0.12559893324407342)
            (6,0.10413381026259408)
            (7,0.08029811114002691)
            (8,0.059179415314921616)
            (9,0.04129197633979064)
            (10,0.03531230963048652)
            (11,0.027030935073055155)
            (12,0.023235192102940394)
            (13,0.02063249259769951)
            (14,0.016286014014564188)
            (15,0.011668456935474683)
            (16,0.00828170107000843)
        }
        ;
    \addplot+[style={solid, very thick}, mesh, mark={none}, point meta={log10(1e-1)}]
        coordinates {

            (2,0.32876010315819804)
            (3,0.1463790642133359)
            (4,0.12881116673624962)
            (5,0.12411814972758203)
            (6,0.10384204572694918)
            (7,0.08114684135625042)
            (8,0.05975051579259385)
            (9,0.04126053401315516)
            (10,0.03469824459246694)
            (11,0.0267222317411153)
            (12,0.022384569393262374)
            (13,0.019973017229876638)
            (14,0.017151604726942712)
            (15,0.013704040669762253)
            (16,0.010938827611209276)
            (17,0.006242067679632966)
        }
        ;
    \addplot+[style={solid, very thick}, mesh, mark={none}, point meta={log10(1e-2)}]
        coordinates {

            (2,0.3243343629354033)
            (3,0.1374243697758083)
            (4,0.12101025312581631)
            (5,0.11929228921903225)
            (6,0.10174685591903625)
            (7,0.08317273087590442)
            (8,0.05964173554568237)
            (9,0.04413182736507329)
            (10,0.0326300360702475)
            (11,0.029034700555228853)
            (12,0.02548981948622611)
            (13,0.021849430828396525)
            (14,0.02026808546161283)
            (15,0.019533956051956236)
            (16,0.018996598691726724)
            (17,0.016719522365180484)
            (18,0.014885885399927815)
            (19,0.013849606422625253)
            (20,0.012666095967953562)
            (21,0.010950296341748067)
            (22,0.011098630109120848)
            (23,0.010968411528122422)
            (24,0.010988908954867814)
            (25,0.011018955774910572)
            (26,0.010941530354057307)
            (27,0.01079222719398965)
            (28,0.010623901633709286)
            (29,0.010368961715245545)
            (30,0.010104620222709904)
            (31,0.009799746978614254)
       }
        ;
    \addplot+[style={solid, very thick}, mesh, mark={none}, point meta={log10(1e-3)}]
        coordinates {

            (2,0.3038177380768124)
            (3,0.16498606501928795)
            (4,0.10334494509234723)
            (5,0.09552643365798621)
            (6,0.09921890824824826)
            (7,0.09584563646250654)
            (8,0.0782471663641827)
            (9,0.06313851945502755)
            (10,0.05224803115795036)
            (11,0.04644028664866159)
            (12,0.04429403959378283)
            (13,0.04428091867292427)
            (14,0.04321762558880941)
            (15,0.041958352053688766)
            (16,0.04127915577677544)
            (17,0.04069187165720609)
            (18,0.039921614999181586)
            (19,0.039394257285715994)
            (20,0.03691946327007171)
            (21,0.035016845608142055)
            (22,0.032230626383386064)
            (23,0.030210235713587528)
            (24,0.02777625404065464)
            (25,0.027286099837435112)
            (26,0.02680483733375886)
            (27,0.02505209234637943)
            (28,0.02280243343592578)
            (29,0.021374017509829583)
            (30,0.020672246836046973)
            (31,0.020227568585349645)
            (32,0.01997796738671822)
            (33,0.019830710124758597)
            (34,0.019547006777054874)
            (35,0.019244061467227482)
            (36,0.018840099809669292)
            (37,0.018380453026738068)
            (38,0.017875525755960582)
            (39,0.01735837536753268)
            (40,0.016971244121199415)
            (41,0.016617411807256458)
            (42,0.016344664181624023)
            (43,0.016103801829297045)
            (44,0.01587411085607682)
            (45,0.015659624197536706)
            (46,0.015431034936358582)
            (47,0.015323242729321012)
            (48,0.01529685500531699)
            (49,0.015246566059793127)
            (50,0.015218342028156528)
            (51,0.015171822644639435)
            (52,0.01510712154747967)
            (53,0.015029239816434265)
            (54,0.014969138072736546)
            (55,0.014985631224805172)
            (56,0.01502326925787081)
            (57,0.015052608986239533)
            (58,0.015045496319238238)
            (59,0.014961486266390942)
            (60,0.014834122626919253)
            (61,0.01457250851274028)
            (62,0.014259467676292022)
            (63,0.013954288129564197)
            (64,0.013656996044476677)
            (65,0.01338218510970453)
            (66,0.013111809592948392)
            (67,0.01286246926777719)
            (68,0.012609878385040593)
            (69,0.01235252512895079)
            (70,0.01215803829564473)
            (71,0.011917003681567326)
            (72,0.011812576303331355)
            (73,0.011801043508235519)
            (74,0.011818645055889473)
            (75,0.011816357190904962)
            (76,0.011669788675910585)
            (77,0.011616813563254336)
            (78,0.01144723934416611)
            (79,0.01148154781560623)
            (80,0.011530221036591128)
            (81,0.011409139794367263)
            (82,0.011301349820173175)
            (83,0.011264662288341585)
            (84,0.011244538574689385)
            (85,0.011212587005153799)
            (86,0.011180554745361776)
            (87,0.011122777342479285)
            (88,0.011034209312752374)
            (89,0.010812669622640709)
            (90,0.010640724980445014)
            (91,0.01053973504315133)
            (92,0.010434546848042382)
            (93,0.010323155427425219)
            (94,0.0101751923789138)
            (95,0.00993150326140518)
        }
        ;
    \end{axis}
\end{tikzpicture}}
    \hfill
    \subfloat[\centering \label{fig:residual_l}\texttt{BWFLnet}, Two-level, $\delta=10 \, \si{\meter}$]{\input{bwfl_residuals_two_level_c}}
    \vspace{0.25cm}
    \caption{Convergence plots for all numerical experiments. Subplots correspond to experiments with different case network, distributed algorithm, and pressure tolerance $\delta$. Colour gradients delineate convergence behaviour for a range of initial penalty parameters $\beta^1$ and the dotted black line represents convergence tolerance $\epsilon_p = 10^{-2}$.}
    \label{fig:convergence}
\end{figure}

\Cref{fig:convergence} shows convergence plots for all numerical experiments using the standard and two-level ADMM algorithms. Each subplot corresponds to an experiment with different case network, distributed algorithm, and pressure range tolerance $\delta$.
Specifically, \Cref{fig:convergence} illustrates the relationship between the average primal residual $(n_n n_t)^{-\frac{1}{2}} \cdot \|r\|$, where $r = h - \bar{h}$, and the number of ADMM iterations. Note that the cumulative number of inner ADMM iterations is used for the two-level algorithm. Within each subplot, initial penalty parameter $\beta^1$ values are differentiated by a discrete colour bar and the dotted black line denotes the primal residual tolerance $\epsilon_p=10^{-2}$, indicating the algorithm's successful termination.

The results show distinct convergence behaviour across the range of $\beta^1$ values. For instance, larger $\beta^1$ values have relatively fast and stable convergence. On the other hand, the number of iterations using small $\beta^1$ values either increased significantly to meet the termination criterion or, in some cases, resulted in oscillatory behaviour and thus a failure to converge within the maximum number of iterations $k_{\max}=10^3$. The latter situation only appeared in experiments involving the \texttt{Modena} network and using the standard ADMM method with $\beta^1 < 10^{-1}$. In these instances, the oscillations' volatility and deviation from the termination criterion intensified as the pressure tolerance decreased. In contrast, the two-level algorithm consistently achieved convergence across all experiments. In exchange for algorithm robustness, however, the two-level algorithm exhibited a relatively slower convergence rate compared to the standard ADMM implementation. This behavior was more pronounced with smaller $\beta^1$ values, especially in the \texttt{Modena} network. \texttt{Modena}'s network characteristics, including its highly looped structure and relatively equidistant PCV actuators, might be resulting in numerous local optima with similar objective values. In contrast, \texttt{BWFLnet} has a relatively branched structure, limiting the number of available flow paths for optimizing pressure or self-cleaning operations. Since water networks vary in topology and source conditions, tuning the penalty parameter for ADMM-based distributed algorithms is crucial for achieving good performance.

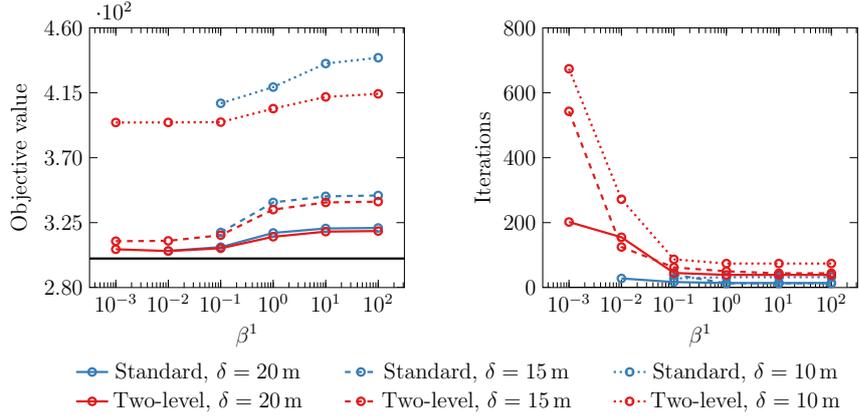
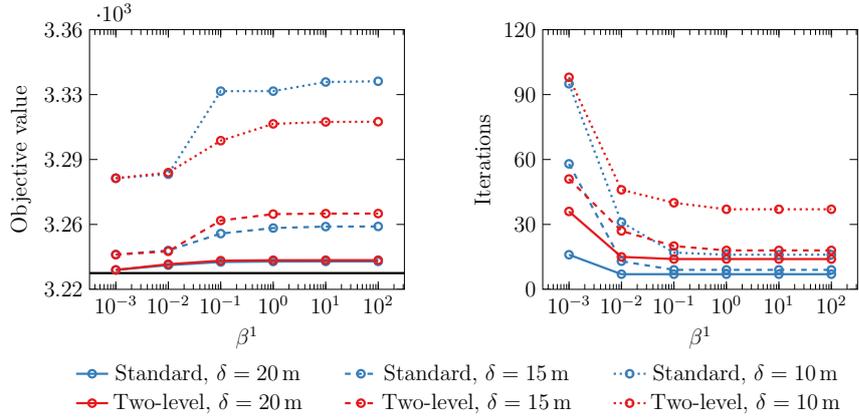
\begin{figure}[t]
    \centering
    \subfloat[\label{fig:obj_iter_modena}\texttt{Modena} ]{% Recommended preamble:
% \pgfplotsset{cycle list/Set1-9}
\begin{tikzpicture}[scale=0.7]
\begin{groupplot}[group style={group size={2 by 1}, horizontal sep={2.6cm}}, cycle list/Set1-3, legend pos=outer north east, legend style={font=\large, row sep=2pt, at={(1.175,-0.26)}, anchor=north, draw=none, /tikz/column 2/.style={column sep=20pt}, /tikz/column 4/.style={column sep=20pt}, cells={align=right}}, legend columns=3, width=7.5cm, height=6.5cm]
    
    \nextgroupplot[xmode={log}, xlabel={$\beta^1$}, ylabel={Objective value}, label style={{font=\large}}, tick label style={{font=\large}}, scaled y ticks={{base 10:-2}}, ymin={280}, ymax={460}, ytick={{280, 325, ..., 460}}, tick style={black}, y tick label style={{/pgf/number format/fixed zerofill}}]
        \addplot+[style={solid, very thick}, mark={o}, mark options={solid}, index of colormap=1 of Set1-3]
        coordinates {

            (0.010000000000000002,305.45077055240256)
            (0.1,308.0069622947647)
            (1.0,317.86241563647775)
            (10.0,320.96278680287924)
            (100.0,321.30484121410024)
        }
        ;
        \addlegendentry {Standard, $\delta=20 \, \si{\meter}$}
        \addplot+[style={dashed, very thick}, mark={o}, mark options={solid}, index of colormap=1 of Set1-3]
        coordinates {

            (0.1,318.17322489497076)
            (1.0,339.03745505063637)
            (10.0,343.23140011411783)
            (100.0,343.82873238628895)
        }
        ;
        \addlegendentry {Standard, $\delta=15 \, \si{\meter}$}
        \addplot+[style={dotted, very thick}, mark={o}, mark options={solid}, index of colormap=1 of Set1-3]
        coordinates {

            (0.1,407.733256104469)
            (1.0,418.9819885597461)
            (10.0,435.2891330963972)
            (100.0,439.3160221555751)
        }
        ;
        \addlegendentry {Standard, $\delta=10 \, \si{\meter}$}
        \addplot+[style={solid, very thick}, mark={o}, mark options={solid}, index of colormap=0 of Set1-3]
            coordinates {
                (0.001,306.53731449055743)
                (0.010000000000000002,305.3003934427188)
                (0.1,307.2301645935022)
                (1.0,315.1972266741197)
                (10.0,318.75359695724813)
                (100.0,319.1720427591751)
            }
            ;
        \addlegendentry {Two-level, $\delta=20 \, \si{\meter}$}
        \addplot+[style={dashed, very thick}, mark={o}, mark options={solid}, index of colormap=0 of Set1-3]
            coordinates {
                (0.001,312.14131138984567)
                (0.010000000000000002,312.37451198222794)
                (0.1,316.292520830823)
                (1.0,334.04566131093475)
                (10.0,338.9947917384118)
                (100.0,339.5468187063409)
            }
            ;
        \addlegendentry {Two-level, $\delta=15 \, \si{\meter}$}
        \addplot+[style={dotted, very thick}, mark={o}, mark options={solid}, index of colormap=0 of Set1-3]
            coordinates {
                (0.001,394.43189715731916)
                (0.010000000000000002,394.46549346342897)
                (0.1,394.6825326485997)
                (1.0,404.04636115364013)
                (10.0,412.16540112247003)
                (100.0,414.275086280649)
            }
            ;
        \addlegendentry {Two-level, $\delta=10 \, \si{\meter}$}
        \draw[style={solid, very thick}, color={black}] ({rel axis cs:1,0}|-{axis cs:0,300.13381858166815}) -- ({rel axis cs:0,0}|-{axis cs:0,300.13381858166815});
        
    \nextgroupplot[xmode={log}, xlabel={$\beta^1$}, ylabel={Iterations}, ymin={0}, ymax={800}, tick style={black}, ytick={{0, 200, ..., 800}}, label style={{font=\large}}, tick label style={{font=\large}}]
        \addplot+[style={solid, very thick}, mark={o}, mark options={solid}, index of colormap=1 of Set1-3, forget plot]
        coordinates {

            (0.010000000000000002,28)
            (0.1,17)
            (1.0,14)
            (10.0,14)
            (100.0,14)
        }
        ;
        \addplot+[style={dashed, very thick}, mark={o}, mark options={solid}, index of colormap=1 of Set1-3, forget plot]
        coordinates {

            (0.1,40)
            (1.0,13)
            (10.0,13)
            (100.0,13)
        }
        ;
        \addplot+[style={dotted, very thick}, mark={o}, mark options={solid}, index of colormap=1 of Set1-3, forget plot]
        coordinates {

            (0.1,28)
            (1.0,32)
            (10.0,32)
            (100.0,32)
        }
        ;
        \addplot+[style={solid, very thick}, mark={o}, mark options={solid}, index of colormap=0 of Set1-3, forget plot]
            coordinates {
                (0.001,202)
                (0.010000000000000002,155)
                (0.1,45)
                (1.0,39)
                (10.0,39)
                (100.0,39)
            }
            ; 
        \addplot+[style={dashed, very thick}, mark={o}, mark options={solid}, index of colormap=0 of Set1-3, forget plot]
            coordinates {
                (0.001,543)
                (0.010000000000000002,125)
                (0.1,62)
                (1.0,50)
                (10.0,44)
                (100.0,44)
            }
            ;
        \addplot+[style={dotted, very thick}, mark={o}, mark options={solid}, index of colormap=0 of Set1-3, forget plot]
            coordinates {
                (0.001,674)
                (0.010000000000000002,272)
                (0.1,87)
                (1.0,74)
                (10.0,74)
                (100.0,74)
            }
            ;
        \addplot+[style={solid, very thick}, color={black}, forget plot]
            coordinates {
                (0.001,-1)
                (0.010000000000000002,-1)
                (0.1,-1)
                (1.0,-1)
                (10.0,-1)
                (100.0,-1)
            }
            ;
    \end{groupplot}

    \end{tikzpicture}
    } \\
    % \vspace{-0.15cm}
    \subfloat[\label{fig:obj_iter_bwfl}\texttt{BWFLnet} ]{% Recommended preamble:
% \pgfplotsset{cycle list/Set1-9}
\begin{tikzpicture}[scale=0.7]
\begin{groupplot}[group style={group size={2 by 1}, horizontal sep={2.6cm}}, cycle list/Set1-3, legend pos=outer north east, legend style={font=\large, row sep=2pt, at={(1.175,-0.26)}, anchor=north, draw=none, /tikz/column 2/.style={column sep=20pt}, /tikz/column 4/.style={column sep=20pt}, cells={align=right}}, legend columns=3, width=7.5cm, height=6.5cm]

    \nextgroupplot[xmode={log}, xlabel={$\beta^1$}, ylabel={Objective value}, label style={{font=\large}}, tick label style={{font=\large}}, legend style={{font=\large}}, scaled y ticks={{base 10:-3}}, ymin={3220}, ymax={3360}, ytick={{3220, 3255, ..., 3360}}, tick style={black}, y tick label style={{/pgf/number format/fixed zerofill}}]
        \addplot+[style={solid, very thick}, mark={o}, mark options={solid, fill opacity=0.15}, index of colormap=1 of Set1-3]
            coordinates {
                (0.001,3230.4586128738915)
                (0.010000000000000002,3233.029718726793)
                (0.1,3234.7390616972907)
                (1.0,3235.016561961062)
                (10.0,3235.055650775417)
                (100.0,3235.059751665121)
            }
            ;
        \addlegendentry {Standard, $\delta=20 \, \si{\meter}$}
        \addplot+[style={dashed, very thick}, mark={o}, mark options={solid, fill opacity=0.15}, index of colormap=1 of Set1-3]
            coordinates {
                (0.001,3238.7930644829985)
                (0.010000000000000002,3240.9479541321975)
                (0.1,3250.0623316076003)
                (1.0,3253.0501147480395)
                (10.0,3253.8361353533805)
                (100.0,3253.916466878546)
            }
            ;
        \addlegendentry {Standard, $\delta=15 \, \si{\meter}$}
        \addplot+[style={dotted, very thick}, mark={o}, mark options={solid, fill opacity=0.15}, index of colormap=1 of Set1-3]
            coordinates {
                (0.001,3279.963969123283)
                (0.010000000000000002,3282.0769499686235)
                (0.1,3326.869000383073)
                (1.0,3326.869000383073)
                (10.0,3331.787335949281)
                (100.0,3332.18904016442)
            }
            ;
        \addlegendentry {Standard, $\delta=10 \, \si{\meter}$}
        \addplot+[style={solid, very thick}, mark={o}, mark options={solid, fill opacity=0.15}, index of colormap=0 of Set1-3]
            coordinates {
                (0.001,3230.456583740568)
                (0.010000000000000002,3233.553680535126)
                (0.1,3235.4627719334167)
                (1.0,3235.6991770591717)
                (10.0,3235.7212357246435)
                (100.0,3235.7235522311503)
            }
            ;
        \addlegendentry {Two-level, $\delta=20 \, \si{\meter}$}
        \addplot+[style={dashed, very thick}, mark={o}, mark options={solid, fill opacity=0.15}, index of colormap=0 of Set1-3]
            coordinates {
                (0.001,3238.674218238799)
                (0.010000000000000002,3240.5082418524653)
                (0.1,3257.102158268221)
                (1.0,3260.5699401781467)
                (10.0,3260.8338687251153)
                (100.0,3260.8621217029468)
            }
            ;
        \addlegendentry {Two-level, $\delta=15 \, \si{\meter}$}
        \addplot+[style={dotted, very thick}, mark={o}, mark options={solid, fill opacity=0.15}, index of colormap=0 of Set1-3]
            coordinates {
                (0.001,3279.872040630907)
                (0.010000000000000002,3282.932330940424)
                (0.1,3300.1545157753444)
                (1.0,3309.2213023201725)
                (10.0,3310.264717354739)
                (100.0,3310.3713441756413)
            }
            ;
        \addlegendentry {Two-level, $\delta=10 \, \si{\meter}$}
        \draw[style={solid, very thick}, color={black}] ({rel axis cs:1,0}|-{axis cs:0,3228.7512850808685}) -- ({rel axis cs:0,0}|-{axis cs:0,3228.7512850808685});
        
        \nextgroupplot[xmode={log}, xlabel={$\beta^1$}, ylabel={Iterations}, ymin={0}, ymax={120}, tick style={black}, ytick={{0, 30, ..., 120}}, label style={{font=\large}}, tick label style={{font=\large}}]
        \addplot+[style={solid, very thick}, mark={o}, mark options={solid, fill opacity=0.15}, index of colormap=1 of Set1-3, forget plot]
            coordinates {
                (0.001,16)
                (0.010000000000000002,7)
                (0.1,7)
                (1.0,7)
                (10.0,7)
                (100.0,7)
            }
            ;
        \addplot+[style={dashed, very thick}, mark={o}, mark options={solid, fill opacity=0.15}, index of colormap=1 of Set1-3, forget plot]
            coordinates {
                (0.001,58)
                (0.010000000000000002,13)
                (0.1,9)
                (1.0,9)
                (10.0,9)
                (100.0,9)
            }
            ;
        \addplot+[style={dotted, very thick}, mark={o}, mark options={solid, fill opacity=0.15}, index of colormap=1 of Set1-3, forget plot]
            coordinates {
                (0.001,95)
                (0.010000000000000002,31)
                (0.1,17)
                (1.0,16)
                (10.0,16)
                (100.0,16)
            }
            ;
        \addplot+[style={solid, very thick}, mark={o}, mark options={solid, fill opacity=0.15}, index of colormap=0 of Set1-3, forget plot]
            coordinates {
                (0.001,36)
                (0.010000000000000002,15)
                (0.1,14)
                (1.0,14)
                (10.0,14)
                (100.0,14)
            }
            ; 
        \addplot+[style={dashed, very thick}, mark={o}, mark options={solid, fill opacity=0.15}, index of colormap=0 of Set1-3, forget plot]
            coordinates {
                (0.001,51)
                (0.010000000000000002,27)
                (0.1,20)
                (1.0,18)
                (10.0,18)
                (100.0,18)
            }
            ;
        \addplot+[style={dotted, very thick}, mark={o}, mark options={solid, fill opacity=0.15}, index of colormap=0 of Set1-3, forget plot]
            coordinates {
                (0.001,98)
                (0.010000000000000002,46)
                (0.1,40)
                (1.0,37)
                (10.0,37)
                (100.0,37)
            }
            ;
        \addplot+[style={solid, very thick}, mark={none}, color={black}]
            coordinates {
                (0.001,-1)
                (0.010000000000000002,-1)
                (0.1,-1)
                (1.0,-1)
                (10.0,-1)
                (100.0,-1)
            }
            ;
    \end{groupplot}
    
\end{tikzpicture}
    } 
    \vspace{0.25cm}
    \caption{Distributed optimization results across a range of initial penalty parameters $\beta^1$. Objective values and number of iterations are shown in the left and right plots, respectively. Pressure tolerances $\delta_i, \, \forall i \in \{20, 15, 10\}$ (in meters) are differentiated by line style. The solid black line is the unconstrained case (i.e. $\bar{\mathcal{X}}=\mathbb{R}$).}
    \label{fig:obj_iter_results}
\end{figure}

In \Cref{fig:obj_iter_results}, we explore the trade-off between objective value performance and iteration count for $\beta^1$ values ranging from $10^{-3}$ to $10^2$.
Inflection points are observed at $\beta^1=10^{-1}$ for \texttt{Modena} and $\beta^1=10^{-2}$ for \texttt{BWFLnet}, marking the best trade-off between objective value and computational efficiency across the tested pressure tolerances. The results also show a strong agreement between the distributed algorithms when the standard ADMM method successfully converges, suggesting that cases with convergence issues can be discarded without compromising solution quality. From an operational perspective, \Cref{fig:obj_iter_results} also highlights a notable reduction in objective value (and marginal decrease in iteration count) corresponding to higher pressure tolerances. The unconstrained case (solid black line) serves as a baseline to assess this improvement. Overall, tuning the penalty parameter is essential for optimizing performance across different networks and problem conditions.

\begin{figure}[h]
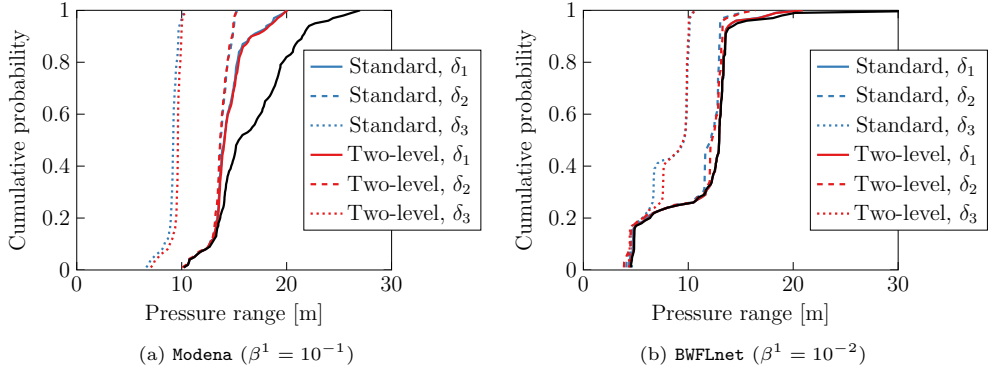

    \centering
    \subfloat[\label{fig:cdf_a}\texttt{Modena} ($\beta^1=10^{-1}$) ]{\input{modena_cdf}}
    \hfill
    \subfloat[\label{fig:cdf_b}\texttt{BWFLnet} ($\beta^1=10^{-2}$)]{\input{bwfl_cdf.tex}} 
    \vspace{0.25cm}
    \caption{Cumulative distribution of temporal pressure range across network nodes. Distributed algorithm is depicted by plot colour and pressure tolerances $\delta_i, \, \forall i \in \{20, 15, 10\}$ (in meters) are differentiated by line style. The solid black line is the unconstrained case (i.e. $\bar{\mathcal{X}}=\mathbb{R}$).}
    \label{fig:cdf_plots}
\end{figure}

The main advantage of the distributed optimization algorithms explored in this work is their fast convergence. It is important, however, to ensure that solution quality is not compromised to achieve this. In \Cref{tab:results_comparison}, we compare solutions computed by the centralized IPOPT solver with those obtained from the standard and two-level ADMM distributed algorithms. Note that the distributed algorithms apply the best $\beta^1$ value derived from our penalty parameter tuning in \Cref{fig:obj_iter_results} and the standard ADMM approach uses a fixed penalty parameter $\rho=2\beta^1$.
\begin{sidewaystable}
\centering
\caption{Comparison of centralized and distributed solution methods. Tolerance $\delta=\infty$ corresponds to the unconstrained pressure range case (i.e. $\bar{\mathcal{X}}=\mathbb{R}$). Dashes denote experiments which failed to converge to a feasible solution within a 1-hour time limit.}
\setlength{\tabcolsep}{12pt}
{\renewcommand{\arraystretch}{1.2}
\small{
\begin{tabular}{lcccclcccc}
    \toprule
    \multicolumn{3}{c}{Experiment} & \multicolumn{2}{c}{Centralized IPOPT} & \multicolumn{5}{c}{Distributed optimization} \\ 
    \cmidrule(lr){1-3} \cmidrule(lr){4-5} \cmidrule(lr){6-10}
    Network & $\beta^1$ & $\delta$ [m] & Obj. & Time [s] & Method & Iter. & Obj. & $\delta_{\text{viol}}$ [m] & Time [s] \\
    \midrule
    \texttt{Modena} & $10^{-1}$ & $\infty$ & $298.3$ & $28.4$ & Standard & $1$ & $300.1$ & $0.00$ & $13.6$ \\ [3pt]
    & & & & & Two-level & $1$ & $300.1$ & $0.00$ & $25.5$ \\ [3pt]
    & & $20$ & $307.7$ & $71.9$ & Standard & $17$ & $308.0$ & $0.00$ & $29.1$ \\ [3pt]
    & & & & & Two-level & $45$ & $307.2$ & $0.06$ & $39.2$ \\ [3pt]
    & & $15$ & $312.9$ & $98.5$ & Standard & $40$ & $318.2$ & $0.32$ & $29.8$ \\ [3pt]
    & & & & & Two-level & $62$ & $316.3$ & $0.32$ & $45.9$ \\ [3pt]
    & & $10$ & $390.8$ & $73.1$ & Standard & $28$ & $407.7$ & $0.21$ & $23.7$ \\ [3pt]
    & & & & & Two-level & $87$ & $394.7$ & $0.43$ & $57.1$ \\ [3pt]
    % \midrule
    \texttt{BWFLnet} & $10^{-2}$& $\infty$ & $3228$ & $1320$ & Standard & $1$ & $3229$ & $0.00$ & $78.0$ \\ [3pt]
    & & & & & Two-level & $1$ & $3229$ & $0.00$ & $79.8$ \\ [3pt]
    & & $20$ & $-$ & $-$ & Standard & $7$ & $3233$ & $0.32$ & $165$ \\ [3pt]
    & & & & & Two-level & $15$ & $3234$ & $0.92$ & $425$ \\ [3pt]
    & & $15$ & $-$ & $-$ & Standard & $13$ & $3241$ & $0.85$ & $291$ \\ [3pt]
    & & & & & Two-level & $27$ & $3241$ & $1.41$ & $763$ \\ [3pt]
    & & $10$ & $-$ & $-$ & Standard & $31$ & $3282$ & $0.71$ & $689$ \\ [3pt]
    & & & & & Two-level & $46$ & $3283$ & $0.58$ & $1210$ \\ [3pt]
    \bottomrule
\end{tabular}
}
}
\label{tab:results_comparison}
\end{sidewaystable}
In addition to objective value and computational time, \Cref{tab:results_comparison} reports the cumulative number of inner ADMM iterations for the two-level algorithm, consistent with the presentation in \Cref{fig:convergence}. Additionally, we tabulate the maximum pressure constraint violations $\delta_{\text{viol}}$. These constraint violations were observed since both distributed algorithms employ a stopping criterion based on the Euclidean norm, which considers the overall distance of the solution from the feasible region rather than the maximum violation, leading to minor constraint violations. More specifically, less than $0.01\%$ of nodes failed to satisfy constraints $\bar{\mathcal{X}}$, with a maximum violation of $0.32 \, \si{\meter}$ and $1.41 \, \si{\meter}$ recorded for \texttt{Modena} and \texttt{BWFLnet}, respectively. It is worth noting that the two-level algorithm yielded slightly larger constraint violations in most experiments. In any case, these violations were considered acceptable due to their order of magnitude aligning with the inherent modeling uncertainties in water networks. Cumulative distribution plots in \Cref{fig:cdf_plots} provide a visualization of the temporal pressure range across the network for different tolerances. These plots highlight the limited impact of the maximum constraint violations within the broader context of the entire network. 

\Cref{tab:results_comparison} highlights the significant computational advantages of the distributed methods, particularly in generating feasible solutions for the large-scale \texttt{BWFLnet} experiments. Among these methods, the two-level algorithm showed only marginal improvements in solution quality over the standard ADMM approach. However, its computational times were typically two to three times longer than those of standard ADMM, especially in experiments involving \texttt{BWFLnet} with smaller pressure tolerances. This increase in computational time is due to the greater iteration complexity required by the two-level algorithm to ensure convergence for nonconvex problems. The two-level algorithm could offer an advantage in networks where frequent re-tuning of the penalty parameters might be necessary for convergence. Despite the longer computational times, the maximum reported time of $20$ minutes across all experiments indicates that the two-level algorithm is a suitable (and robust) algorithm for near real-time (e.g. hourly) control in large-scale water networks.

%%% Conclusion %%%
\section{Conclusion}
\label{sec:conclusion}
We presented a new control model for optimizing pressure and water quality operations in water distribution networks. The formulation imposed a set of time-coupling constraints to manage temporal pressure variations, which are exacerbated by the transition between pressure and water quality control objectives. We then investigated distributed optimization methods to solve the resulting time-coupled, nonconvex problem, as state-of-the-art nonlinear solvers struggled to find feasible solutions. In particular, we implemented two distributed optimization algorithms based on the alternating direction method of multipliers (ADMM): a standard ADMM scheme and a two-level variant that provided convergence guarantees for our nonconvex problem.

We evaluated the distributed algorithms using a benchmarking water network and a large-scale operational network in the UK. Our numerical experiments first demonstrated the sensitivity of convergence behaviour for different choices of the initial penalty parameter. This was highlighted by the standard ADMM method failing to converge for a number of experiments. Although the two-level algorithm proved to be more robust, finding feasible solutions to all numerical experiments, computational times increased significantly for small initial penalty parameters. With an appropriately tuned penalty parameter, however, both distributed algorithms yielded good quality solutions and computational times compatible with near real-time (e.g. hourly) control implementations. In future work, we aim to embed such efficient and robust solution methods within a model predictive control framework to optimize water network operations in near real-time.

%%==================================%%
%%             POSTSCRIPT           %%
%%==================================%%

\backmatter

\section*{Declarations}
\bmhead{Funding}
Bradley Jenks is supported by Imperial College London (Skempton Scholarship), Bristol Water Plc, Analytical Technology, and the Natural Sciences and Engineering Research Council of Canada (PGSD-577767-2023). Ivan Stoianov is supported by a Royal Academy of Engineering Senior Research Fellowship in Dynamically Adaptive Water Supply Networks (RCSRF2324-17-41).

\bmhead{Data availability}
Data and code associated with this manuscript can be accessed at the following GitHub repository: \url{https://github.com/bradleywjenks/wdn_control_admm.git}.

\bmhead{Author contributions}
All authors contributed to the study conception and design. Optimization and analysis scripts were prepared by Bradley Jenks. The initial draft of the manuscript was written by Bradley Jenks and all authors provided feedback on previous versions. The final manuscript was reviewed and approved by all authors.

\bmhead{Competing interests}
The authors declare no competing interests.

\begin{appendices}

\section{Two-level convergence properties}\label{sec:A1}

The following presents theoretical convergence guarantees for solving problem \eqref{eq:problem_reform_2} to an approximate stationary solution using the two-level algorithm.
\vspace{9pt}

\begin{remark}
\label{math:remark_1}
We show that Assumptions 1--3 provided in \citet{SUN2023} hold for problem \eqref{eq:problem_reform_2}:
\begin{itemize}
    \item \textbf{Assumption 1.} Problem \eqref{eq:problem_reform_2} is feasible given that the chosen pressure range tolerance $\delta$ in \eqref{eq:pr_constraint} can be satisfied. We ensure $\delta$ is large enough by using the natural pressure variation from diurnal loading conditions as a baseline (i.e. no control).
    \item \textbf{Assumption 2.} The objective function $f$, which is a sum of linear and nonlinear functions $f_{\text{AZP}}$ and $f_{\text{SCC}}$, respectively, is continuously differentiable. Moreover, the sets $\mathcal{X}_t$ and $\bar{\mathcal{X}}$ are compact and defined by a finite number of constraints ${p_t = 3n_p+2n_n+2n_v+2n_f}$ and ${\bar{p} = n_tn_n + n_n}$, respectively. The condition that $\bar{\mathcal{X}}$ is convex is also met since it is comprised of linear constraints \eqref{eq:pr_constraint_linear}.
    \item \textbf{Assumption 3.} The first ADMM block update in \eqref{eq:3block_admm_x_update} can find a stationary solution via efficient local optimization methods \citep[Section 1.4.1.]{BOYD2004}. In this work, we use state-of-the-art nonlinear solver IPOPT \citep{WACHTER2006}, which applies an interior point line search filter method to compute local solutions to the nonlinear subproblems. Since we initialize ADMM with a feasible starting point (see line 5 in \Cref{alg:two_level_algorithm}), the updated iterate $(q^{k+1}, h^{k+1}, u^{k+1}, \bar{h}^k, z^k, y^k)$ is a stationary point of \eqref{eq:3block_admm_x_update} with improved objective value. Thus, 
    \begin{equation}
        \hat{L}^m(q^{k+1},h^{k+1},u^{k+1},\bar{h}^k,z^k,y^k) \leq \hat{L}^m(q^k,h^k,u^k,\bar{h}^k,z^k,y^k) < \infty
    \end{equation}
    holds at outer level iteration $m$ for all $k \in \mathbb{Z}_{++}$.
\end{itemize}

\end{remark}

\noindent
Under Assumptions 1--3 in \Cref{math:remark_1}, we can establish convergence guarantees for the two-level algorithm following \citep[Proposition 1 and Theorem 2]{SUN2023}. Specifically, if $\rho^m=2\beta^m$, then the sequence of solutions $\{(q^{m}, h^{m},u^{m},\bar{h}^m,z^{m},y^{m})\}_{m \in \mathbb{Z}_{++}}$ converge to a limit point $(q^{*}, h^{*},u^{*},\bar{h}^{*},z^{*},y^{*})$.

\newpage
\section{Nomenclature}

\begin{table}[h!]
\centering
\captionsetup{width=1\linewidth}
\caption{List of key symbols and definitions.}
\setlength{\tabcolsep}{12pt}
{\renewcommand{\arraystretch}{1.2}
\small{
\begin{tabular}{cl}
\toprule
Symbol          & Definition \\ \midrule
$\mathcal{N}$   & Set of model nodes\\
$\mathcal{P}$   & Set of model links\\
$\mathcal{T}$   & Set of time steps in control horizon\\
$n_n, n_0$      & Number of junction and source nodes\\
$n_f$           & Number of AFV nodes\\
$n_p$           & Number of links\\
$n_v$           & Number of PCV links\\
$n_t$           & Number of discrete time steps\\
$A_{10}, A_{12}$& Link-source and -junction node incidence matrices\\
$A_{13}, A_{14}$& Matrix mappings of PCV and AFV locations\\
$d_t$           & Nodal demands at time $t$\\
$h_{0t}$        & Source hydraulic heads at time $t$\\
$q_t, h_t$      & Flow and hydraulic head variables at time $t$\\
$u_t = (\eta_t, \alpha_t)$ & PCV and AFV control variables at time $t$\\
$h^{\min}$     & Minimum regulatory pressure head\\
$\alpha^{\max}$ & Maximum flow rate at AFV nodes\\
$f_t(\cdot)$    & Control problem objective function at time $t$\\
$\zeta$         & Node elevations\\
$\psi(\cdot)$   & Sigmoid function in approximated SCC objective\\
$\mathcal{X}_t$ & Set of hydraulic constraints at time $t$\\
$\mathcal{\bar{X}}$ & Set of time-coupling constraints\\
$\delta$        & Pressure range tolerance\\
$\bar{h}$       & Duplicated hydraulic head variables for consensus constraint \eqref{eq:global_copy}\\
$z$             & Slack variables for ALR problem reformulation\\
$L(\cdot)$      & Augmented Lagrangian function\\
$y$             & Dual variables associated with constraints \eqref{eq:prob_reform_1b} and \eqref{eq:prob_alr_b}\\
$\lambda$       & Dual variables associated with constraint \eqref{eq:prob_reform_2d}\\
$\rho$          & Fixed ADMM penalty parameter\\
$\beta$         & Adaptive ALM penalty parameter\\
$\gamma, \omega$& $\beta$ update parameters\\
$\epsilon$      & Numerical convergence tolerance\\

\bottomrule
\end{tabular}
}
}
\label{tab:centralized_solver}
\end{table}

\end{appendices}

\newpage

\bibliography{references}

\begin{thebibliography}{42}
\providecommand{\natexlab}[1]{#1}
\providecommand{\url}[1]{{#1}}
\providecommand{\urlprefix}{URL }
\providecommand{\doi}[1]{\url{https://doi.org/#1}}
\providecommand{\eprint}[2][]{\url{#2}}
 \bibcommenthead

\bibitem[{Abraham et~al(2016)Abraham, Blokker, and Stoianov}]{ABRAHAM2016}
Abraham E, Blokker M, Stoianov I (2016) Network analysis, control valve placement and optimal control of flow velocity for self-cleaning water distribution systems. 18th Conference on Water Distribution System Analysis, WDSA2016 pp 1--9. \doi{10.1016/j.proeng.2017.03.272}

\bibitem[{Abraham et~al(2018)Abraham, Blokker, and Stoianov}]{ABRAHAM2018}
Abraham E, Blokker M, Stoianov I (2018) Decreasing the discoloration risk of drinking water distribution systems through optimized topological changes and optimal flow velocity control. Journal of Water Resources Planning and Management 144(2):04017093. \doi{10.1061/(asce)wr.1943-5452.0000878}

\bibitem[{Bezanson et~al(2017)Bezanson, Edelman, Karpinski, and Shah}]{BEZANSON2017}
Bezanson J, Edelman A, Karpinski S, et~al (2017) Julia: A fresh approach to numerical computing. SIAM Review 59(1):65--98. \doi{10.1137/141000671}

\bibitem[{Blokker et~al(2010)Blokker, Vreeburg, Schaap, and van Dijk}]{BLOKKER2010}
Blokker M, Vreeburg J, Schaap P, et~al (2010) The self-cleaning velocity in practice. Water Distribution System Analysis (WDSA) pp 187--199. \doi{10.1061/41203(425)19}

\bibitem[{Boxall et~al(2023)Boxall, Blokker, Schaap, Speight, and Husband}]{BOXALL2023}
Boxall J, Blokker M, Schaap P, et~al (2023) Managing discolouration in drinking water distribution systems by integrating understanding of material behaviour. Water Research 243:120416. \doi{10.1016/j.watres.2023.120416}

\bibitem[{Boyd and Vandenberghe(2004)}]{BOYD2004}
Boyd S, Vandenberghe L (2004) Convex optimization. Cambridge University Press

\bibitem[{Boyd et~al(2010)Boyd, Parikh, Chu, Peleato, and Eckstein}]{BOYD2010}
Boyd S, Parikh N, Chu E, et~al (2010) Distributed optimization and statistical learning via the alternating direction method of multipliers. Foundations and Trends in Machine Learning 3(1):1--122. \doi{10.1561/2200000016}

\bibitem[{Bragalli et~al(2012)Bragalli, Lodi, and D'Ambrosio}]{BRAGALLI2012}
Bragalli C, Lodi A, D'Ambrosio C (2012) On the optimal design of water distribution networks: a practical {MINLP} approach. Optimization and Engineering 13(2):219--246. \doi{10.1007/s11081-011-9141-7}

\bibitem[{Bui et~al(2022)Bui, Jeong, and Kang}]{BUI2022}
Bui XK, Jeong G, Kang D (2022) Adaptive {DMA} design and operation under multiscenarios in water distribution networks. Sustainability 14(6):3692. \doi{10.3390/su14063692}

\bibitem[{Dunning et~al(2017)Dunning, Huchette, and Lubin}]{DUNNING2017}
Dunning I, Huchette J, Lubin M (2017) {JuMP}: A modeling language for mathematical optimization. SIAM Review 59(2):295--320. \doi{10.1137/15M1020575}

\bibitem[{Eckstein and Bertsekas(1992)}]{ECKSTEIN1992}
Eckstein J, Bertsekas DP (1992) On the {Douglas-Rachford} splitting method and the proximal point algorithm for maximal monotone operators. Mathematical Programming 55:293--318. \doi{10.1007/BF01581204}

\bibitem[{Fooladivanda and Taylor(2018)}]{FOOLADIVANDA2018}
Fooladivanda D, Taylor JA (2018) Energy-optimal pump scheduling and water flow. IEEE Transactions on Control of Network Systems 5(3):1016--1026. \doi{10.1109/TCNS.2017.2670501}

\bibitem[{Ghaddar et~al(2015)Ghaddar, Naoum-Sawaya, Kishimoto, Taheri, and Eck}]{GHADDAR2015}
Ghaddar B, Naoum-Sawaya J, Kishimoto A, et~al (2015) A lagrangian decomposition approach for the pump scheduling problem in water networks. European Journal of Operational Research 241(2):490--501. \doi{10.1016/j.ejor.2014.08.033}

\bibitem[{Gholami et~al(2023)Gholami, Sun, Zhang, and Sun}]{GHOLAMI2023}
Gholami A, Sun K, Zhang S, et~al (2023) An {ADMM-based} distributed optimization method for solving security-constrained {AC} optimal power flow. INFORMS Operations Research \doi{10.1287/opre.2023.2486}

\bibitem[{Giudicianni et~al(2020)Giudicianni, Herrera, di~Nardo, and Adeyeye}]{GIUDICIANNI2020}
Giudicianni C, Herrera M, di~Nardo A, et~al (2020) Automatic multiscale approach for water networks partitioning into dynamic district metered areas. Water Resources Management 34(2):835--848. \doi{10.1007/s11269-019-02471-w}

\bibitem[{{Gurobi Optimization}(2023)}]{GUROBI2023}
{Gurobi Optimization} (2023) {Gurobi Optimizer Reference Manual, Version 10.0.2}

\bibitem[{Hong et~al(2016)Hong, Luo, and Razaviyayn}]{HONG2016}
Hong M, Luo ZQ, Razaviyayn M (2016) Convergence analysis of alternating direction method of multipliers for a family of nonconvex problems. SIAM Journal on Optimization 26(1):337--364. \doi{10.1137/140990309}

\bibitem[{HSL(2021)}]{HSL2021}
HSL (2021) A collection of fortran codes for large scale scientific computation. \urlprefix\url{http://www.hsl.rl.ac.uk/}

\bibitem[{Jara-Arriagada and Stoianov(2021)}]{ARRIAGADA2021}
Jara-Arriagada C, Stoianov I (2021) Pipe breaks and estimating the impact of pressure control in water supply networks. Reliability Engineering and System Safety 210:107525. \doi{10.1016/j.ress.2021.107525}

\bibitem[{Jenks et~al(2023{\natexlab{a}})Jenks, Pecci, and Stoianov}]{JENKS2023a}
Jenks B, Pecci F, Stoianov I (2023{\natexlab{a}}) Optimal design-for-control of self-cleaning water distribution networks using a convex multi-start algorithm. Water Research 231:119602. \doi{10.1016/j.watres.2023.119602}

\bibitem[{Jenks et~al(2023{\natexlab{b}})Jenks, Ulusoy, Pecci, and Stoianov}]{JENKS2023b}
Jenks B, Ulusoy AJ, Pecci F, et~al (2023{\natexlab{b}}) Dynamically adaptive networks for integrating optimal pressure management and self-cleaning controls. Annual Reviews in Control 55:486--497. \doi{10.1016/j.arcontrol.2023.03.014}

\bibitem[{Jiang et~al(2019)Jiang, Lin, Ma, and Zhang}]{JIANG2019}
Jiang B, Lin T, Ma S, et~al (2019) Structured nonconvex and nonsmooth optimization: algorithms and iteration complexity analysis. Computational Optimization and Applications 72(1):115--157. \doi{10.1007/s10589-018-0034-y}

\bibitem[{Li and Pong(2015)}]{LI2015}
Li G, Pong TK (2015) Global convergence of splitting methods for nonconvex composite optimization. SIAM Journal on Optimization 25(4):2434--2460. \doi{10.1137/140998135}

\bibitem[{Magnusson et~al(2016)Magnusson, Weeraddana, Rabbat, and Fischione}]{MAGNUSSON2016}
Magnusson S, Weeraddana PC, Rabbat MG, et~al (2016) On the convergence of alternating direction lagrangian methods for nonconvex structured optimization problems. IEEE Transactions on Control of Network Systems 3(3):296--309. \doi{10.1109/TCNS.2015.2476198}

\bibitem[{Martínez-Codina et~al(2016)Martínez-Codina, Castillo, González-Zeas, and Garrote}]{MARTINEZ2016}
Martínez-Codina, Castillo M, González-Zeas D, et~al (2016) Pressure as a predictor of occurrence of pipe breaks in water distribution networks. Urban Water Journal 13(7):676--686. \doi{10.1080/1573062X.2015.1024687}

\bibitem[{Nitivattananon et~al(1996)Nitivattananon, Sadowski, and Quimpo}]{NITIV1996}
Nitivattananon V, Sadowski EC, Quimpo RG (1996) Optimization of water supply system operation. Journal of Water Resources Planning and Management 122(5):374--384. \doi{10.1061/(ASCE)0733-9496(1996)122:5(374)}

\bibitem[{Pas et~al(2022)Pas, Schuurmans, and Patrinos}]{PAS2022}
Pas P, Schuurmans M, Patrinos P (2022) {ALPAQA}: A matrix-free solver for nonlinear {MPC} and large-scale nonconvex optimization. 2022 European Control Conference (ECC) pp 417--422. \doi{10.23919/ECC55457.2022.9838172}

\bibitem[{Pecci et~al(2017)Pecci, Abraham, and Stoianov}]{PECCI2017b}
Pecci F, Abraham E, Stoianov I (2017) Scalable pareto set generation for multiobjective co-design problems in water distribution networks: a continuous relaxation approach. Structural and Multidisciplinary Optimization 55(3):857--869. \doi{10.1007/s00158-016-1537-8}

\bibitem[{Rezaei et~al(2015)Rezaei, Ryan, and Stoianov}]{REZAEI2015}
Rezaei H, Ryan B, Stoianov I (2015) Pipe failure analysis and impact of dynamic hydraulic conditions in water supply networks. Procedia Engineering 119(1):253--262. \doi{10.1016/j.proeng.2015.08.883}

\bibitem[{Schwaller and van Zyl(2015)}]{SCHWALLER2015}
Schwaller J, van Zyl JE (2015) Modeling the pressure-leakage response of water distribution systems based on individual leak behavior. Journal of Hydraulic Engineering 141(5). \doi{10.1061/(asce)hy.1943-7900.0000984}

\bibitem[{Sun and Sun(2021)}]{SUN2021}
Sun K, Sun X (2021) A two-level {ADMM algorithm} for {AC OPF} with global convergence guarantees. IEEE Transactions on Power Systems 36(6):5271--5281. \doi{10.1109/TPWRS.2021.3073116}

\bibitem[{Sun and Sun(2023)}]{SUN2023}
Sun K, Sun X (2023) A two-level distributed algorithm for nonconvex constrained optimization. Computational Optimization and Applications 84(2):609--649. \doi{10.1007/s10589-022-00433-4}

\bibitem[{Tang and Daoutidis(2022)}]{TANG2022}
Tang W, Daoutidis P (2022) Fast and stable nonconvex constrained distributed optimization: the {ELLADA} algorithm. Optimization and Engineering 23(1):259--301. \doi{10.1007/s11081-020-09585-w}

\bibitem[{Themelis and Patrinos(2020)}]{THEMELIS2018}
Themelis A, Patrinos P (2020) {Douglas–Rachford} splitting and {ADMM} for nonconvex optimization: Tight convergence results. SIAM Journal on Optimization 30(1):149--181. \doi{10.1137/18M1163993}

\bibitem[{Ulusoy et~al(2023)Ulusoy, Nerantzis, and Stoianov}]{ULUSOY2023}
Ulusoy AJ, Nerantzis D, Stoianov I (2023) Adaptive {MPC} for burst incident management in water distribution networks. IEEE Transactions on Control of Network Systems \doi{10.1109/TCNS.2023.3259103}

\bibitem[{Vreeburg et~al(2009)Vreeburg, Blokker, Horst, and van Dijk}]{VREEBURG2009}
Vreeburg J, Blokker M, Horst P, et~al (2009) Velocity-based self-cleaning residential drinking water distribution systems. Water Science and Technology: Water Supply 9(6):635--641. \doi{10.2166/ws.2009.689}

\bibitem[{Wang et~al(2019)Wang, Yin, and Zeng}]{WANG2019}
Wang Y, Yin W, Zeng J (2019) Global convergence of {ADMM} in nonconvex nonsmooth optimization. Journal of Scientific Computing 78(1):29--63. \doi{10.1007/s10915-018-0757-z}

\bibitem[{Wright et~al(2014)Wright, Stoianov, Parpas, Henderson, and King}]{WRIGHT2014}
Wright R, Stoianov I, Parpas P, et~al (2014) Adaptive water distribution networks with dynamically reconfigurable topology. Journal of Hydroinformatics 16(6):1280--1301. \doi{10.2166/hydro.2014.086}

\bibitem[{Wright et~al(2015)Wright, Abraham, Parpas, and Stoianov}]{WRIGHT2015}
Wright R, Abraham E, Parpas P, et~al (2015) Control of water distribution networks with dynamic {DMA} topology using strictly feasible sequential convex programming. Water Resources Research 51(12):9925--9941. \doi{10.1002/2015WR017466}

\bibitem[{Wächter and Biegler(2006)}]{WACHTER2006}
Wächter A, Biegler LT (2006) On the implementation of an interior-point filter line-search algorithm for large-scale nonlinear programming. Mathematical Programming 106(1):25--57. \doi{10.1007/s10107-004-0559-y}

\bibitem[{Zamzam et~al(2019)Zamzam, Dall'Anese, Zhao, Taylor, and Sidiropoulos}]{ZAMZAM2019}
Zamzam AS, Dall'Anese E, Zhao C, et~al (2019) Optimal water-power flow-problem: Formulation and distributed optimal solution. IEEE Transactions on Control of Network Systems 6(1):37--47. \doi{10.1109/TCNS.2018.2792699}

\bibitem[{Zessler and Shamir(1989)}]{ZESSLER1989}
Zessler BU, Shamir U (1989) Optimal operation of water distribution systems. Journal of Water Resources Planning and Management 115(6):735--752. \doi{10.1061/(ASCE)0733-9496(1989)115:6(735)}

\end{thebibliography}

\end{document}